\documentclass[12pt]{article}

\usepackage{indentfirst}
\usepackage{amsfonts}
\usepackage{amsmath}
\usepackage{amssymb}
\usepackage{amsbsy, amsthm}
\usepackage[margin=1in]{geometry}
\usepackage{imakeidx,xcolor}

\usepackage{graphicx}
\usepackage{caption}
\usepackage{subcaption}

\newtheorem{dfn}{Definition} [section]
\newtheorem{theorem}[dfn]{Theorem}
\newtheorem{lemma}[dfn]{Lemma}

\newtheorem{corollary}[dfn]{Corollary}

\newenvironment{pf}{\noindent{\bf Proof.}}
{\enspace\vrule height5pt depth0pt width5pt}
\def\deg {{\rm deg}}
\def\ch {{\rm ch}}
\def\leng {{\rm leng}}

\title{Odd list-coloring of graphs of small Euler genus with no short cycles of specific types}
\author{Rishi Balaji\thanks{mail2mlrb2@gmail.com. Stanford Online High School, Redwood City, CA.}, Victoria Khazhinsky\thanks{victoria.khazhinsky@gmail.com. LASA High School, Austin, TX.}, Chun-Hung Liu\thanks{chliu@tamu.edu. Department of Mathematics, Texas A\&M University, USA. Partially supported by NSF under CAREER award DMS-2144042.}, Kevin Qin\thanks{kevinqinhshl@gmail.com. College Station High School, TX, USA.}}

\begin{document}

\maketitle

\begin{abstract}
Odd coloring is a variant of proper coloring and has received wide attention. 
We study the list-coloring version of this notion in this paper.
We prove that if $G$ is a graph embeddable in the torus or the Klein bottle with no cycle of length 3, 4, and 6 such that no 5-cycles share an edge, then for every function $L$ that assigns each vertex of $G$ a set $L(v)$ of size 5, there exists a proper coloring that assigns each vertex $v$ of $G$ an element of $L(v)$ such that for every non-isolated vertex, some color appears an odd number of times on its neighborhood.
In particular, every graph embeddable in the torus or the Klein bottle with no cycle of length 3, 4, 6, and 8 is odd 5-choosable.
The number of colors in these results are optimal, and there exist graphs embeddable in those surfaces of girth 6 requiring six or seven colors.
\end{abstract}

\section{Introduction}

A {\it proper coloring} of a graph\footnote{Graphs are finite and simple in this paper. In other words, every graph has a finite number of vertices but does not contain a loop or parallel edges.} is a function that assigns each vertex a color such that no two adjacent vertices receive the same color.
The {\it chromatic number} $\chi(G)$ of a graph $G$ is the minimum $k$ such that $G$ admits a proper coloring of $k$ colors.
Upper bounds for the chromatic number of graphs in terms of their Euler genus have been well-studied.
They can often be reduced if we forbid certain lengths of cycles.
For example, planar graphs without 3-cycle are properly 3-colorable \cite{g}, compared to the general upper bound 4 provided by the Four Color Theorem \cite{ah1,ahk,ah2,rst}; similarly, graphs of Euler genus $g$ with no 3-cycle are properly $O((g/\log g)^{1/3})$-colorable \cite{gt}, improving the general upper bound $O(\sqrt{g})$ \cite{ah,d,h1,r}.
Moreover, Erd\H{o}s suggested to study the minimum $k$ such that planar graphs with no cycle of length between $4$ and $k$ are 3-colorable as an approach to Steinberg's conjecture (see \cite{s}).
The currently best result is $6 \leq k \leq 7$, where the lower bound was proved in \cite{chkls} and disproved Steinberg's conjecture, and the upper bound was proved in \cite{bgrs}.
Related results, such as requiring cycles of certain lengths that satisfy certain adjacency or distance conditions rather than forbidding them, were also studied \cite{br,x}.

In this paper, we consider odd coloring of graphs embedded in surfaces with additional constraints on their cycles. 
Odd coloring of graphs with certain topological properties has attracted wide attention \cite{cps,cckp,c,cls,dmo,h,lwy,m,nz,pp,ps,ty,wy} although this concept was only introduced fairly recently in \cite{ps}.
See \cite{alo,accchhhkz,cckp2,cckp3,cl,dop,kp,l,lr,lr2} for related results.

Let $G$ be a graph and $k$ a positive integer. 
A {\it $k$-coloring} of $G$ is a function $f: V(G) \rightarrow \{1,2,...,k\}$.
A $k$-coloring is {\it odd} if it is proper, and for every non-isolated vertex $v$, some color appears an odd number of times in the neighborhood of $v$.
We say that $G$ is {\it odd $k$-colorable} if $G$ admits an odd $k$-coloring.

It is easy to see that the 5-cycle is a planar graph that is odd 5-colorable but not odd 4-colorable, since the neighbors of every degree-2 vertex must receive different colors in every odd coloring.
Petru\v{s}evski and \v{S}krekovski \cite{ps} conjectured that every planar graph is odd 5-colorable and proved that every planar graph is odd 9-colorable.
This result was improved in two directions: Petr and Portier \cite{pp} proved that every planar graph is odd 8-colorable, and Metrebian \cite{m} and Tian and Yin \cite{ty} proved that every graph embeddable in the torus is odd $9$-colorable.

Note that, unlike proper coloring, degeneracy and bipartiteness are not sufficient to ensure an odd coloring with a small number of colors.
Such examples can be obtained by considering the 1-subdivision of another graph.
We say that a graph $G$ is a {\it 1-subdivision} of $H$ if $G$ is obtained from $H$ by subdividing every edge exactly once.
If $G$ is a 1-subdivision of $H$, then the restriction of every odd coloring of $G$ on $V(H)$ is a proper coloring of $H$, so $G$ is not odd-$(\chi(H)-1)$-colorable.
By choosing $H$ to be graphs of arbitrarily large chromatic number, we know that the resulting 1-subdivisions of $H$ require arbitrarily large numbers of colors for odd coloring, even though the 1-subdivision of $H$ is 2-degenerate and bipartite no matter what $H$ is.

On the other hand, degeneracy still provides an upper bound for the required number of colors for odd coloring if we restrict ourselves to minor-closed families.
For example, a corollary of a result in \cite{l} implies that graphs in any $d$-degenerate minor-closed family are odd $(2d+1)$-colorable. 
For every fixed surface $\Sigma$, the graphs embeddable in $\Sigma$ form a minor-closed family, and Euler's formula provides an upper bound for the degeneracy of graphs embedded in $\Sigma$ and hence provide an upper bound for the number of required colors for odd coloring. 
In addition, Cranston \cite{c} proved that every graph with maximum average degree less than $\frac{4c}{c+2}$ is odd $c$-colorable for every $c \in \{5,6\}$, which implies that every planar graph of girth at least 7 is odd 5-colorable; 
Cho, Choi, Kwon, and Park \cite{cckp} complemented this result by proving the cases $c \geq 7$. 

In this paper, we not only consider graphs embeddable in a surface but also consider list-coloring.
Let $G$ be a graph.
For every positive integer $k$, a {\it $k$-list-assignment} $L$ of $G$ is a collection $(L_v: v \in V(G))$ of sets with $|L_v|=k$ for every $v \in V(G)$.
For a $k$-list-assignment $L$ of $G$, an {\it $L$-coloring} is a function $f$ that maps each vertex $v$ of $G$ to an element of $L_v$; an $L$-coloring is {\it proper} if $u$ and $v$ receive different colors for every edge $uv$ of $G$.
We say that $G$ is {\it odd $k$-choosable} if for every $k$-list-assignment $L$ of $G$, there exists a proper $L$-coloring such that for every non-isolated vertex $v$, some color appears an odd number of times in the neighborhood of $v$.

Clearly, every odd $k$-choosable graph is odd $k$-colorable.
But the converse is not true even for bipartite graphs with girth at least 8.
For example, it is not hard to see that for every integer $d$, $K_{d,d}$ is not $k$-choosable for some $k=\Omega(\log d)$; so the 1-subdivision of $K_{d,d}$ is not odd $k$-choosable for some $k=\Omega(\log d)$, even though it is easy to show that it is odd $4$-colorable and has girth at least 8.

The following is the main result of this paper.

\begin{theorem} \label{main_simple}
If $G$ is a graph embeddable in the torus or the Klein bottle with no cycle of length 3, 4, and 6 such that no two 5-cycles share an edge, then $G$ is odd 5-choosable.
\end{theorem}

Note that the number of colors mentioned in Theorem \ref{main_simple} is optimal since every 5-cycle is not odd 4-colorable.
Moreover, the 1-subdivision of $K_7$ (and $K_6$, respectively) is not odd 6-colorable (and 5-colorable, respectively), but it has girth 6 and can be embedded in the torus (and the Klein bottle, respectively).
So the condition that there is no 6-cycle in Theorem \ref{main_simple} cannot be dropped. 
In addition, Theorem \ref{main_simple} strengthens the aforementioned corollary of a result of Cranston \cite{c} that every planar graph of girth at least 7 is odd 5-colorable. 

In addition, if there exist two 5-cycles sharing an edge, then there exists a cycle of length 3, 4, 6, or 8.
So Theorem \ref{main_simple} immediately leads to the following corollary.

\begin{corollary}
If $G$ is a graph embeddable in the torus or the Klein bottle with no cycle of length 3, 4, 6, and 8, then $G$ is odd 5-choosable.
\end{corollary}

We in fact prove a stronger result in order to prove Theorem \ref{main_simple}.

Let $G$ be a graph and $R \subseteq E(G)$.
We define the {\it $R$-length} of a cycle $C$ in $G$ to be $|E(C)|+|E(C) \cap R|$.
We say a vertex $v$ of $G$ is {\it $R$-relaxed} if either the degree of $v$ is odd or $0$, or $v$ is incident with an edge in $R$. 

\begin{theorem} \label{main_real}
Let $G$ be a graph embeddable in the torus or the Klein bottle, and let $R$ be a subset of $E(G)$.
If no cycle of $G$ has $R$-length equal to 3, 4, or 6, and no two cycles of $R$-length 5 share exactly one edge, then for every 5-list-assignment $L$ of $G$, there exists a proper $L$-coloring of $G$ such that for every non-$R$-relaxed vertex $v$, some color appears an odd number of times on the neighborhood of $v$. 
\end{theorem}

Clearly, for any vertex $v$ of odd degree, some color appears an odd number of times on its neighborhood in any proper coloring.
So Theorem \ref{main_simple} follows from Theorem \ref{main_real} by taking $R=\emptyset$.

We will prove Theorem \ref{main_real} by using the discharging method.
We first introduce basic terminology in Section \ref{sec:terminology}.
Then we study structures of a minimal counterexample in Section \ref{sec:reducible} and apply discharging arguments in Section \ref{sec:discharging} to derive a contradiction.

\section{Terminology} \label{sec:terminology}

Let $G$ be a graph, and let $v$ be a vertex of $G$.
We denote the set of neighbors of $v$ by $N_G(v)$, and denote the degree of $v$ by $\deg_G(v)$.
We define $N_G[v]$ to be $N_G(v) \cup \{v\}$.
If the graph $G$ is clear from the context, we simply denote $N_G(v)$, $\deg_G(v)$, and $N_G[v]$ by $N(v)$, $\deg(v)$, and $N[v]$, respectively.
Recall that a coloring $f$ is {\it odd} if for every vertex $v$ with $N(v) \neq \emptyset$, some color appears an odd number of times on $N(v)$.

Let $S$ be a subset of $V(G)$.
We denote the subgraph of $G$ induced by $S$ by $G[S]$.
For every list-assignment $L$ of $G$, we define $L|_S$ to be the list-assignment $(L_v: v \in S)$ of $G[S]$.

Let $k$ be an integer.
A {\it $k$-vertex} is a vertex of degree $k$.
A {\it $(\geq k)$-vertex} is a vertex of degree at least $k$.
Now we assume that $G$ is embedded in a surface.
We define the {\it length} of a face $f$, denoted by $\leng(f)$, to be the minimum number of edges of a closed walk that contains all edges incident with $f$.
A {\it $k$-face} is a face of length $k$.
A {\it $(\geq k)$-face} is a face of length at least $k$.
We say that two faces {\it share} edges if some edge is incident with both faces; we say that two faces are {\it adjacent} if they share edges.

\section{Structures of a minimal counterexample} \label{sec:reducible}

In this section, we assume that $G$ is a counterexample to Theorem \ref{main_real} such that $|V(G)|$ is as small as possible.
Clearly, $|V(G)| \geq 2$.

\begin{lemma} \label{conn}
$G$ is connected.
\end{lemma}

\begin{pf}
If $G$ is not connected, then each component has a proper $L$-coloring such that for every non-$R$-relaxed vertex $v$, some color appears an odd number of times on $N_G(v)$.
By combining such colorings for all components, we obtain a proper $L$-coloring of $G$ such that for every non-$R$-relaxed vertex $v$, some color appears an odd number of times on $N_G(v)$.
\end{pf}

\bigskip

Since $G$ is connected (by Lemma \ref{conn}) and embeddable in the torus or the Klein bottle, $G$ has a 2-cell embedding in a surface of Euler genus at most $2$. 
We fix such a 2-cell embedding of $G$ in this section.

\begin{lemma} \label{min_deg}
Every vertex of $G$ has degree at least three.
\end{lemma}

\begin{pf}
Let $v$ be a vertex of minimum degree in $G$.
It suffices to show that $\deg(v) \geq 3$.
Suppose to the contrary that $\deg(v) \leq 2$.
By Lemma \ref{conn}, $v$ has degree at least 1 since $|V(G)| \geq 2$.

\medskip

\noindent{\bf Claim 1:} $\deg(v) \geq 2$.

\noindent{\bf Proof of Claim 1:}
Suppose that $\deg(v)=1$.
Let $u$ be the unique neighbor of $v$.
If $\deg(u)=1$, then $G=K_2$, so it is clear that $G$ has a proper $L$-coloring of $G$ such that every non-$R$-relaxed vertex $x$, some color appears an odd number of times on $N_{G}(x)$, a contradiction.

Hence $\deg(u) \geq 2$.
By the minimality of $|V(G)|$, $G-v$ has a proper $L|_{V(G-v)}$-coloring $f_1$ such that for every non-$(R \cap E(G-v))$-relaxed vertex $x$ in $G-v$, some color $c_x$ appears an odd number of times on $N_{G-v}(x)$.
Note that for every non-$R$-relaxed vertex $x$ in $G$, if $x \not \in \{u,v\}$, then since $N_G(x)=N_{G-v}(x)$, we know that $x$ is non-$(R \cap E(G-v))$-relaxed in $G-v$.
Since $\deg(u) \geq 2$, if $u$ is non-$R$-relaxed in $G$, then either $u$ is non-$(R \cap E(G-v))$-relaxed in $G-v$ or $\deg_{G-v}(u)$ is odd.
If $u$ is $R$-relaxed in $G$, then let $c_u$ be an arbitrary color; if $\deg_{G-v}(u)$ is odd, then some $c_u$ appears an odd number of times on $N_{G-v}(u)$.
By further coloring $v$ with a color in $L_v-\{f_1(u),c_u\}$, we obtain a proper $L$-coloring of $G$ such that every non-$R$-relaxed vertex $x$, some color appears an odd number of times on $N_{G}(x)$, a contradiction.
$\Box$

\medskip

So $\deg(v)=2$.
In particular, every vertex of $G$ has degree at least two.
Let $N_G(v)=\{a,b\}$.

\medskip

\noindent{\bf Claim 2:} $v$ is not incident with an edge in $R$.

\noindent{\bf Proof of Claim 2:}
Suppose that $v$ is incident with an edge in $R$.
By symmetry, we may assume that $av \in R$.
By the minimality of $|V(G)|$, $G-v$ has a proper $L|_{V(G-v)}$-coloring $f_2$ such that for every non-$(R \cap E(G-v))$-relaxed vertex $x$ in $G-v$, some color $c_x$ appears an odd number of times on $N_{G-v}(x)$.
If $b$ is $R$-relaxed, then let $c_b$ be an arbitrary color.
If $b$ is non-$R$-relaxed and $(R \cap E(G-v))$-relaxed, then $\deg_{G-v}(b)$ is odd, so some color $c_b$ appears an odd number of times on $N_{G-v}(x)$.

Since $av \in R$, we know that $a$ and $v$ are $R$-relaxed.
So by further coloring $v$ with a color in $L_v-\{f_2(a),f_2(b),c_b\}$, we obtain a proper $L$-coloring of $G$ such that every non-$R$-relaxed vertex $x$, some color appears an odd number of times on $N_{G}(x)$, a contradiction.
$\Box$

\medskip

Hence $v$ is not incident with an edge in $R$.
If $ab \in R$, then $abv$ is a cycle in $G$ of $R$-length 4, a contradiction.
If $ab \in E(G)-R$, then $abv$ is a cycle in $G$ of $R$-length 3, a contradiction.
So $ab \not \in E(G)$.

Let $G'$ be the graph obtained from $G-v$ by adding an edge $ab$, and let $R' = R \cup \{ab\}$.

\medskip

\noindent{\bf Claim 3:} No cycle in $G'$ has $R'$-length equal to 3, 4, or 6, and no two cycles in $G'$ of $R'$-length 5 share exactly one edge in $G'$.

\noindent{\bf Proof of Claim 3:}
For every cycle $C'$ of $G'$, if $C'$ does not contain $ab$, then $C'$ is also a cycle in $G$ with $R'$-length equal to its $R$-length; if $C'$ contains $ab$, then $(C'-ab)+avb$ is a cycle of $G$ with $R$-length equal to the $R'$-length of $C'$.
So no cycle of $G'$ has $R'$-length equal to 3, 4, or 6.

Suppose that there exist two cycles $C_1'$ and $C_2'$ in $G'$ of $R'$-length 5 sharing exactly one edge in $G'$.
If $E(C_1')$ and $E(C_2')$ do not share $ab$, then there exist two cycles in $G$ of $R$-length 5 sharing exactly one edge in $G$, a contradiction.
So the edge shared by $C_1'$ and $C_2'$ is $ab$.
Since $ab \in R'$, the cycle $(C_1'-ab) \cup (C_2'-ab)$ is a cycle in $G$ of $R$-length $(5-2)+(5-2)=6$, a contradiction.
$\Box$

\medskip

Since $|V(G')|<|V(G)|$, the minimality of $|V(G)|$ implies that there exists a proper $L|_{V(G-v)}$-coloring $f_3$ such that for every non-$R'$-relaxed vertex $x$ in $G'$, some color appears an odd number of times on $N_{G'}(x)$.
For every $x \in \{a,b\}$, if some color appears an odd number of times on $N_{G'-(\{a,b\}-\{x\})}(x)$, then let $\ell_x$ be such a color; otherwise, let $\ell_x$ be an arbitrary color.
Now we further color $v$ with a color in $L_v-\{f_3(a),f_3(b),\ell_x,\ell_y\}$ to obtain an $L$-coloring $f^*$ of $G$.
Clearly, $f^*$ is a proper $L$-coloring.
Since $G$ is a counterexample, there exists a non-$R$-relaxed vertex $z$ in $G$ such that no color appears an odd number of times on $N_G(z)$.
Since $ab \in E(G')$, we know $f_3(a) \neq f_3(b)$, so $z \neq v$.
For every $x \in V(G)-\{a,b,v\}$, we know $N_G(x)=N_{G'}(x)$.
So $z \in \{a,b\}$.
If $\deg_G(z)$ is odd, then some color must appear an odd number of times on $N_G(z)$.
So $\deg_G(z)$ is even.
Hence $|N_{G'-(\{a,b\}-\{z\})}|$ is odd, so $\ell_z$ appears an odd number of times on $N_{G'-(\{a,b\}-\{z\})}$.
Since $f^*(v) \neq \ell_z$, we know that $\ell_z$ appears an odd number of times on $N_G(z)$, a contradiction.
This proves the lemma.
\end{pf}

\begin{lemma} \label{3_vertex_relaxed}
No 3-vertex is adjacent to at least two $R$-relaxed vertices.
\end{lemma}

\begin{pf}
Suppose to the contrary that there exists a 3-vertex $v$ adjacent to at least two $R$-relaxed vertices $x$ and $y$.
Let $z$ be the vertex in $N(v)-\{x,y\}$.
Since $|V(G-v)|<|V(G)|$, there exists a proper $L|_{V(G-v)}$-coloring $f$ of $G-v$ such that for every non-$(R \cap E(G-v))$-relaxed vertex $x$ in $G-v$, some color $c_x$ appears an odd number of times on $N_{G-v}(x)$.
If $z$ is $(R \cap E(G-v))$-relaxed in $G$ and $\deg_{G-v}(z)$ is odd, then some color $c_z$ appears an odd number of times on $N_{G-v}(z)$; if $z$ is $(R \cap E(G-v))$-relaxed in $G$ and $\deg_{G-v}(z)$ is even, then let $c_z$ be an arbitrary color.
Now we further color $v$ with a color in $L_v-\{f(x),f(y),f(z),c_z\}$ to obtain an $L$-coloring $f^*$ of $G$.
Clearly, $f^*$ is proper.

Since $G$ is a counterexample, there exists a non-$R$-relaxed vertex $w$ such that no color appears an odd number of times on $N_G(w)$.
Since $\deg_G(v)=3$, $w \neq v$.
If $w \in V(G)-N_G(v)$, then $N_G(w)=N_{G-v}(w)$, so $c_w$ appears an odd number of times on $N_G(w)$.
Hence $w \in N_G(v)$.
Then $w=z$ since $x$ and $y$ are $R$-relaxed.
Since $z$ is non-$R$-relaxed, we know that $\deg_G(z)$ is even and $\deg_{G-v}(z)$ is odd.
Since $f^*(v) \neq c_z$, $c_z$ appears an odd number of times on $N_G(z)$.
So $f^*$ shows that $G$ is not a counterexample, a contradiction.
\end{pf}

\begin{lemma} \label{3_cycle_relaxed}
Every 3-cycle has $R$-length 5, and every vertex contained in a 3-cycle is $R$-relaxed.
\end{lemma}

\begin{pf}
Let $C$ be a 3-cycle.
So the $R$-length of $C$ equals $3+|R \cap E(C)| \not \in \{3,4,6\}$.
Hence $|R \cap E(C)| \not \in \{0,1,3\}$.
So $C$ contains exactly two edges of $R$ and hence has $R$-length 5.
It implies that every vertex of $C$ is incident with an edge in $R$, so it is $R$-relaxed.
\end{pf}

\begin{lemma} \label{3_not_3}
No 3-vertex is contained in a 3-cycle.
\end{lemma}

\begin{pf}
If some 3-vertex is contained in a 3-cycle, then it is adjacent to at least two $R$-relaxed vertices by Lemma \ref{3_cycle_relaxed}, contradicting Lemma \ref{3_vertex_relaxed}.
\end{pf}

\begin{lemma} \label{relaxed_4_neighbors}
No $R$-relaxed 4-vertex is adjacent to four $R$-relaxed vertices.
\end{lemma}

\begin{pf}
Let $v$ be an $R$-relaxed 4-vertex.
Since $G$ is a minimal counterexample, there exists a proper $L|_{V(G-v)}$-coloring $f$ of $G-v$ such that for every non-$(R \cap E(G-v))$-relaxed vertex $x$, some color appears an odd number of times on $N_{G-v}(x)$.
Since all neighbors of $v$ are $R$-relaxed, by coloring $v$ with a color in $L_v-\{f(x): x \in N_G(v)\}$, we obtain a proper $L$-coloring of $G$ such that for every non-$R$-relaxed vertex $x$, some color appears an odd number of times on $N_{G}(x)$, a contradiction.
\end{pf}

\begin{lemma} \label{relaxed_433}
There does not exist a 4-vertex $x$ adjacent to two 3-vertices $y$ and $z$ such that some vertex in $N(y)-\{x\}$ is $R$-relaxed, some vertex in $N(z)-\{x\}$ is $R$-relaxed, and some vertex in $N[x]-\{y,z\}$ is $R$-relaxed. 
\end{lemma}

\begin{pf}
Suppose to the contrary that such vertices $x,y,z$ exist.
Let $G' = G-\{x,y,z\}$.
By the minimality of $G'$, there exists a proper $L|_{V(G')}$-coloring $\phi$ of $G'$ such that for every non-$(R \cap E(G'))$-relaxed vertex, some color appears an odd number of times on its neighborhood in $G'$.
For every vertex $a \in V(G')$, if $a$ is $R$-relaxed, then let $c_a$ be an arbitrary color not in $L_x$; otherwise, let $c_a$ be a color appearing on $N_{G'}(a)$ such that $c_a$ appears an odd number of times on $N_{G'}(a)$ if possible.

Denote $N_G(y)$ by $\{y_1,y_2,x\}$, and denote $N_G(z)$ by $\{z_1,z_2,x\}$.
By assumption, we may assume that $y_1$ and $z_1$ are $R$-relaxed.
Let $S_y = L_y-\{\phi(y_1),\phi(y_2),c_{y_2}\}$, and let $S_z = L_z-\{\phi(z_1),\phi(z_2),c_{z_2}\}$.
Note that $|S_y| \geq 2$ and $|S_z| \geq 2$.

Denote $N_G(x)$ by $\{x_1,x_2,y,z\}$.
Let $S_x = L_x-\{\phi(x_1),\phi(x_2),c_{x_1},c_{x_2}\}$.
By the choice of $c_{x_1}$ and $c_{x_2}$, if $x_1$ or $x_2$ is $R$-relaxed, then $|S_x| \geq 2$.
Since either $x$ is $R$-relaxed, or at least one of $x_1$ and $x_2$ is $R$-relaxed, we know that either $x$ is $R$-relaxed and $|S_x| \geq 1$, or $|S_x| \geq 2$.

\medskip

\noindent{\bf Claim 1:} $\{x_1,x_2\} \cap \{y_1,y_2,z_1,z_2\} = \emptyset$, and $c_a$ appears an odd number of times on $N_{G'}(a)$ for every vertex $a$ in $G'$ non-$R$-relaxed in $G$. 

\noindent{\bf Proof of Claim 1:}
If $x_i=y_j$ for some $i \in [2]$ and $j \in [2]$, then $xx_iyx$ is a triangle containing the 3-vertex $y$, contradicting Lemma \ref{3_not_3}.
Similarly, $x_i \neq z_j$ for any $i \in [2]$ and $j \in [2]$. 
So $\{x_1,x_2\} \cap \{y_1,y_2,z_1,z_2\} = \emptyset$.

Suppose that there exists $v \in V(G')$ non-$R$-relaxed in $G$ such that $c_v$ appears an even number of times on $N_{G'}(v)$. 
By the definition of $c_v$, $v$ is $R \cap E(G')$-relaxed and has even degree in $G'$.
Since $v$ is non-$R$-relaxed in $G$, we know that $v$ is an isolated vertex in $G'$.
So $N_G(v) \subseteq \{x,y,z\}$.
By Lemma \ref{min_deg}, $N_G(v) = \{x,y,z\}$.
So $y$ is a 3-vertex contained in the triangle $xyvx$, contradicting Lemma \ref{3_not_3}.
$\Box$

\medskip

\noindent{\bf Claim 2:} $x$ is non-$R$-relaxed in $G$.

\noindent{\bf Proof of Claim 2:}
Suppose to the contrary that $x$ is $R$-relaxed.
Then we extend $\phi$ to an $L$-coloring $\phi_1$ of $G$ by defining $\phi_1(x)$ to be an element of $S_x$, and for each $a \in \{y,z\}$, defining $\phi_1(a)$ to be an element in $S_a-\{\phi_1(x)\}$.
Clearly, $\phi_1$ is a proper $L$-coloring of $G$.

Since $G$ is a counterexample, there exists a non-$R$-relaxed vertex $v$ such that no color appears an odd number of times on $N_G(v)$.
Since $x$ is $R$-relaxed by assumption, $v \neq x$.
Since 3-vertices are $R$-relaxed, $v \not \in \{x,y,z\}$.
So $v \in V(G')$.
By Claim 1, $c_v$ appears on $N_{G'}(v)$ an odd number of times.
The definition of $\phi$ implies that $c_v$ appears an odd number of times on $N_G(v)$, a contradiction.
$\Box$

\medskip

By Claim 2, $x$ is non-$R$-relaxed in $G$.
Hence $|S_x| \geq 2$.
By the assumption of this lemma, we may assume that $x_1$ is $R$-relaxed in $G$ by symmetry.

\medskip

\noindent{\bf Claim 3:} $\phi(x_1) \neq \phi(x_2)$.

\noindent{\bf Proof of Claim 3:}
Suppose $\phi(x_1)=\phi(x_2)$.
Then $|S_x| \geq 3$.
We extend $\phi$ to an $L$-coloring $\phi_2$ of $G$ by defining $\phi_2(y)$ to be an element of $S_y$, defining $\phi_2(z)$ to be an element of $S_z-\{\phi_2(y)\}$, and defining $\phi_2(x)$ to be an element of $S_x-\{\phi_2(y),\phi_2(z)\}$.
Since $\phi(x_1)=\phi(x_2)$, we know that $\phi_2(y)$ or $\phi_2(z)$ appears exactly once on $N_G(x)$.
So $\phi_2$ shows that $G$ is not a counterexample, a contradiction. 
$\Box$

\medskip

So $\phi(x_1) \neq \phi(x_2)$.

\medskip

\noindent{\bf Claim 4:} $S_y \cap S_z = \emptyset$.

\noindent{\bf Proof of Claim 4:}
Suppose $S_y \cap S_z \neq \emptyset$.
We extend $\phi$ to an $L$-coloring $\phi_3$ of $G$ by defining $\phi_3(y)$ and $\phi_3(z)$ to be an element of $S_y \cap S_z$, and defining $\phi_3(x)$ to be an element of $S_x-\{\phi_3(y),\phi_3(z)\} = S_x-\{\phi_3(y)\} \neq \emptyset$.
Then $\phi_3$ shows that $G$ is not a counterexample since $\phi(x_1)$ or $\phi(x_2)$ appears exactly once on $N_G(x)$, a contradiction.
$\Box$

\medskip

So $S_y \cap S_z = \emptyset$.
Hence $S_y-\{\phi(x_1),\phi(x_2)\} \neq \emptyset$ or $S_z-\{\phi(x_1),\phi(x_2)\} \neq \emptyset$.
By symmetry, we may assume $S_y-\{\phi(x_1),\phi(x_2)\} \neq \emptyset$.
We extend $\phi$ to an $L$-coloring $\phi_4$ of $G$ by defining $\phi_4(y)$ to be an element of $S_y-\{\phi(x_1),\phi(x_2)\}$, defining $\phi_4(x)$ to be an element of $S_x-\{\phi_4(y)\}$, and defining $\phi_4(z)$ to be an element of $S_z-\{\phi_4(x)\}$.

Since $y$ is a 3-vertex, $yz \not \in E(G)$ by Lemma \ref{3_not_3}.
So $\phi_4$ is a proper $L$-coloring of $G$.
Moreover, $\phi_4$ shows that $G$ is not a counterexample since $\phi_4(y)$ appears exactly once on $N_G(x)$, a contradiction.
This proves the lemma.
\end{pf}

\begin{lemma} \label{no_3333_neighbors}
There does not exist a 4-vertex $v$ such that every neighbor of $v$ is a 3-vertex adjacent to an $R$-relaxed vertex other than $v$.
\end{lemma}

\begin{pf}
Suppose to the contrary that there exists a 4-vertex $v$ such that for every $z \in N_G(v)$, there exists $z' \in N_G(z)-\{v\}$ such that $z'$ is $R$-relaxed.
Let $X=\{x \in V(G): N_G(x) \subseteq N_G(v)\}$.
Note that $X \neq \emptyset$ since $v \in X$.
By Lemma \ref{min_deg}, every vertex in $X$ has degree at least three.
By Lemma \ref{3_vertex_relaxed}, no vertex in $X$ is a 3-vertex.
So every vertex in $X$ is a 4-vertex.

Let $G'=G-(X \cup N_G(v))$.
Since $X$ is exactly the set of isolated vertices of $G-N_G(v)$, we know that $G'$ does not have an isolated vertex.

By the minimality of $|V(G)|$, there exists a proper $L|_{V(G')}$-coloring $\phi$ of $G'$ such that for every non-$(R \cap E(G'))$-relaxed vertex $v$ in $G'$, some color appears an odd number of times on $N_{G'}(v)$.
For every $z \in V(G')$, let $c_{z}$ be a color appearing on $N_{G'}(z)$ such that $c_{z}$ appears on $N_{G'}(z)$ an odd number of times if possible; note that if $z$ is non-$R$-relaxed in $G$, then since $G'$ does not have an isolated vertex, either $z$ is non-$(R \cap E(G'))$-relaxed in $G'$, or $\deg_{G'}(z)$ is odd, so $c_{z'}$ appears an odd number of times on $N_{G'}(z')$.

For every $z \in N_G(v)$, let $z''$ be the vertex in $N_G(z)-\{v,z'\}$.
For every $z \in N_G(v)$, let $S_z = L_z - (\{c_{z''}\} \cup \{\phi(u): u \in N_G(z) \cap V(G')\})$; note that $|S_z| \geq 2$.

Let $u \in N_G(v)$.
We extend $\phi$ by further defining $\phi(u)$ to be an element of $S_u$, and then for every $z \in N_G(v)-\{u\}$, defining $\phi(z)$ to be an element of $S_z-\{\phi(u)\}$, and finally for every $x \in X$, defining $\phi(x)$ to be an element of $L_x-\{\phi(z): z \in N_G(x)\}$.
Clearly, $\phi$ is a proper $L$-coloring of $G$.
By the choice of $c_z$ for $z \in V(G')$ and the definition of $S_z$ for $z \in N_G(v)$, we know that $c_z$ appears an odd number of times on $N_G(z)$ for every vertex $z$ in $V(G')$ that is $R$-relaxed in $G$.
Since $N_G(x)=N_G(v)$ for every $x \in X$, we know $\phi(u)$ appears exactly once on $N_G(x)$.
Since every vertex in $N_G(v)$ is a 3-vertex, $\phi$ is an odd $L$-coloring.
\end{pf}

\begin{lemma} \label{relaxed_443333}
There do not exist two adjacent 4-vertices $x$ and $y$ such that for every $z \in \{x,y\}$, $z$ is adjacent to two 3-vertices $z_1$ and $z_2$ such that for every $i \in [2]$, some vertex $z_i'$ in $N(z_i)-\{z\}$ is $R$-relaxed.
\end{lemma}

\begin{pf}
Suppose to the contrary that such 4-vertices $x,y$, 3-vertices $x_1,x_2,y_1,y_2$ and the neighbors $x_1',x_2',y_1',y_2'$ exist.
Note that $x_1,x_2,y_1,y_2$ are pairwise distinct since no 3-vertex is contained in a 3-cycle.

Let $G' = G-\{x,y,x_1,x_2,y_1,y_2\}$.
If there exists an isolated vertex $v$ of $G'$ such that $v$ is non-$R$-relaxed in $G$, then $\deg_G(v)$ is even, so $\deg_G(v) \geq 4$ by Lemma \ref{min_deg}, and hence $N_G(v)=\{x_1,x_2,y_1,y_2\}$ by Lemma \ref{3_not_3}, contradicting Lemma \ref{no_3333_neighbors} (since $v \not \in \{x_1',x_2',y_1',y_2'\}$ as $v$ is non-$R$-relaxed).
So no non-$R$-relaxed vertex in $G$ is an isolated vertex of $G'$.

By the minimality of $|V(G)|$, there exists a proper $L|_{V(G')}$-coloring $\phi$ of $G'$ such that for every non-$(R \cap E(G'))$-relaxed vertex $v$ in $G'$, some color appears an odd number of times on $N_{G'}(v)$.

For every $z \in V(G')$, let $c_{z}$ be a color appearing on $N_{G'}(z)$ such that $c_{z}$ appears on $N_{G'}(z)$ an odd number of times if possible; note that if $z$ is non-$R$-relaxed, then since $z$ is not an isolated vertex of $G'$, either $\deg_{G'}(z)$ is odd or $z$ is non-$(R \cap E(G'))$-relaxed in $G'$, so $c_{z'}$ appears an odd number of times on $N_{G'}(z')$.

Let $x'$ be the vertex in $N(x)-\{y,x_1,x_2\}$, and let $y'$ be the vertex in $N(y)-\{x,y_1,y_2\}$.
For every $z \in \{x,y\}$, let $S_z$ be a subset of $L_z - \{\phi(z'),c_{z'}\}$ of size 3. 
For every $z \in \{x_1,x_2,y_1,y_2\}$, let $z''$ be the vertex in $N_G(z)-\{x,y,z'\}$, and let $S_z$ be a subset of $L_z - (\{c_{z''}\} \cup \{\phi(v): v \in N_G(z) \cap V(G')\})$ of size 2. 

By the definition of $S_z$ for $z \in \{x,y,x_1,x_2,y_1,y_2\}$, for every proper $L$-coloring $\phi'$ of $G$ obtained by extending $\phi$ such that $\phi'(z) \in S_z$ for every $z \in \{x,y,x_1,x_2,y_1,y_2\}$, if there exists a non-$R$-relaxed vertex of $G$ such that no color appears an odd number of times on its neighborhood, then this vertex must be in $\{x,y\}$ since $x_1,x_2,y_1,y_2$ are 3-vertices.

\medskip

\noindent{\bf Claim 1:} $S_{x_1} \cap S_{x_2} \neq \emptyset$ and $S_{y_1} \cap S_{y_2} \neq \emptyset$.

\noindent{\bf Proof of Claim 1:} 
Suppose to the contrary that $S_{x_1} \cap S_{x_2} = \emptyset$.
We extend $\phi$ by defining $\phi(y)$ to be an element of $S_y$, and then for every $z \in \{y_1,y_2\}$, defining $\phi(z)$ to be an element of $S_z-\{\phi(y)\}$.
Since $|S_x|=3$ and $\deg_{G-x}(y)=3$, we can define $\phi(x)$ to be an element of $S_x-\{\phi(y)\}$ such that some color appears on $N_G(y)$ an odd number of times.
Since $S_{x_1} \cap S_{x_2} = \emptyset$, by symmetry we may assume $\phi(x) \not \in S_{x_2}$.
Then we further extend $\phi$ by defining $\phi(x_1)$ to be an element of $S_{x_1}-\{\phi(x)\}$, and defining $\phi(x_2)$ to be an element of $S_{x_2}$ such that some color appears an odd number of times on $N_G(x_2)$.
Clearly, $\phi$ is a proper $L$-coloring.

Since $G$ is a counterexample, there exists a non-$R$-relaxed vertex $v$ in $G$ such that no color appears an odd number of times on $N_G(v)$.
Note that $v \in \{x,y\}$.
But $v \not \in \{x,y\}$ by the definition of $\phi$, a contradiction.

So $S_{x_1} \cap S_{x_2} \neq \emptyset$.
Similarly, $S_{y_1} \cap S_{y_2} \neq \emptyset$.
$\Box$

\medskip

\noindent{\bf Claim 2:} For every $z \in \{x,y\}$, if $\phi'$ is an $L$-coloring of $G$ such that $\phi'(\bar{z}) \neq \phi'(z')$ and $\phi'(z_1)=\phi'(z_2)$, where $\bar{z}$ is the element of $\{x,y\}-\{z\}$, then some color appears an odd number of times on $N_G(z)$.

\noindent{\bf Proof of Claim 2:} 
Since $\phi'(z_1)=\phi'(z_2)$ and $\phi'(\bar{z}) \neq \phi'(z')$, if $\phi'(z_1)=\phi'(z')$, then $\phi'(\bar{z})$ appears on $N_G(z)$ exactly once; if $\phi'(z_1) \neq \phi'(z')$, then $\phi'(z')$ appears on $N_G(z)$ exactly once.
$\Box$

\medskip

\noindent{\bf Claim 3:} $S_{x_1} \neq S_{x_2}$ and $S_{y_1} \neq S_{y_2}$.

\noindent{\bf Proof of Claim 3:} 
Suppose to the contrary that $S_{x_1} = S_{x_2}$.

We extend $\phi$ by defining $\phi(y)$ to be an element of $S_y-\{\phi(x')\}$, and then for every $z \in \{y_1,y_2\}$, defining $\phi(z)$ to be an element of $S_z-\{\phi(y)\}$, and then defining $\phi(x)$ to be an element of $S_x-\{\phi(y)\}$ such that some color appears an odd number of times on $N_G(y)$.
Since $S_{x_1}=S_{x_2}$, we can further define $\phi(x_1)=\phi(x_2)$ to be an element of $S_{x_1}-\{\phi(x)\}$.
Clearly, $\phi$ is a proper $L$-coloring.

Since $G$ is a counterexample, there exists a non-$R$-relaxed vertex $v$ in $G$ such that no color appears an odd number of times on $N_G(v)$.
Note that $v \in \{x,y\}$.
By the definition of $\phi$, $v \neq y$.
By Claim 2, $v \neq x$, a contradiction.

So $S_{x_1} \neq S_{x_2}$.
Similarly, $S_{y_1} \neq S_{y_2}$.
$\Box$

\medskip

\noindent{\bf Claim 4:} $\{\phi(x'),c_{x'}\} \cap S_{x_1} \cap S_{x_2} = \emptyset$.

\noindent{\bf Proof of Claim 4:} 
Suppose to the contrary that there exists $c \in \{\phi(x'),c_{x'}\} \cap S_{x_1} \cap S_{x_2}$.
Since $c \in \{\phi(x'),c_{x'}\}$, we know $c \not \in S_x$.

We extend $\phi$ by defining $\phi(x_1)=\phi(x_2)=c$, and then defining $\phi(y)$ to be an element of $S_y$ such that some color appears an odd number of times on $N_G(x)$, and then for every $z \in \{y_1,y_2\}$, defining $\phi(z)$ to be an element of $S_z-\{\phi(y)\}$.
Since $|S_x| \geq 3$, we can further define $\phi(x)$ to be an element of $S_x-\{\phi(y)\}$ such that some color appears an odd number of times on $N_G(y)$.
Note that $\phi$ is an $L$-coloring since $c \not \in S_x$.

Since $G$ is a counterexample, there exists a non-$R$-relaxed vertex $v$ in $G$ such that no color appears an odd number of times on $N_G(v)$.
Then $v \in \{x,y\}$.
But $v \not \in \{x,y\}$ by the definition of $\phi$.
$\Box$

\medskip

By Claims 1 and 3, there exists a unique element $q$ in $S_{x_1} \cap S_{x_2}$, and there exists a unique element $q'$ in $S_{y_1} \cap S_{y_2}$.

\medskip

\noindent{\bf Claim 5:} $\{\phi(y'),q\} \subseteq S_x$ and $\{\phi(x'),q'\} \subseteq S_y$.

\noindent{\bf Proof of Claim 5:} 
Suppose to the contrary that $\{\phi(y'),q\} \not \subseteq S_x$.

We extend $\phi$ by defining $\phi(y_1)=\phi(y_2)=q'$ and defining $\phi(y)$ to be an element of $S_y-\{\phi(x'),q'\}$. 
Since $\{\phi(y'),q\} \not \subseteq S_x$, we know $|S_x-\{\phi(y'),q\}| \geq |S_x|-1 \geq 2$.
So we can define $\phi(x)$ to be an element of $S_x-\{\phi(y),\phi(y'),q\}$.
Finally, for every $z \in \{x_1,x_2\}$, define $\phi(z)=q$. 
Clearly, $\phi$ is a proper coloring.

Since $G$ is a counterexample, there exists a non-$R$-relaxed vertex $v$ in $G$ such that no color appears an odd number of times on $N_G(v)$.
Note that $v \in \{x,y\}$.
But $v \not \in \{x,y\}$ by Claim 2, a contradiction.

So $\{\phi(y'),q\} \subseteq S_x$.
Similarly, $\{\phi(x'),q'\} \subseteq S_y$.
$\Box$

\medskip

\noindent{\bf Claim 6:} $\phi(y') \neq q$ and $\phi(x') \neq q'$.

\noindent{\bf Proof of Claim 6:} 
Suppose to the contrary that $\phi(y') = q$.

We extend $\phi$ by defining $\phi(x_1)=\phi(x_2)=q$, defining $\phi(y_1)=\phi(y_2)=q'$, defining $\phi(y)$ to be an element of $S_y-\{q',\phi(x')\}$, and defining $\phi(x)$ to be the element of $S_x-\{q,\phi(y)\}$.
Clearly, $\phi$ is a proper $L$-coloring.

Since $G$ is a counterexample, there exists a non-$R$-relaxed vertex $v$ in $G$ such that no color appears an odd number of times on $N_G(v)$.
Then $v \in \{x,y\}$.
But $v \not \in \{x,y\}$ by Claim 2 since $\phi(x) \neq q=\phi(y')$, a contradiction.

So $\phi(y') \neq q$.
Similarly, $\phi(x') \neq q'$.
$\Box$

\medskip

\noindent{\bf Claim 7:} $\phi(x') \in (S_{x_1} \cup S_{x_2}) - \{q\}$.

\noindent{\bf Proof of Claim 7:} 
Suppose to the contrary that $\phi(x') \not \in (S_{x_1} \cup S_{x_2}) - \{q\}$.

By Claim 5, $q \in S_x$.
We extend $\phi$ by defining $\phi(x)=q \in S_x$, and for every $z \in \{x_1,x_2\}$, defining $\phi(z)$ to be the element of $S_z-\{q\}$.
Define $\phi(y_1)=\phi(y_2)=q'$, and define $\phi(y)$ to be the element of $S_y-\{\phi(x),q'\}$.
Clearly, $\phi$ is a proper $L$-coloring.

Since $G$ is a counterexample, there exists a non-$R$-relaxed vertex $v$ in $G$ such that no color appears an odd number of times on $N_G(v)$.
Note that $v \in \{x,y\}$.
By Claim 6, $\phi(x)=q \neq \phi(y')$.
So $v \neq y$ by Claim 2.
Hence $v=x$.
Since $\phi(x'),\phi(x_1),\phi(x_2)$ are pairwise distinct, some color appears exactly once on $N_G(x)$, a contradiction.
$\Box$

\medskip

By Claim 7 and symmetry, we may assume $\phi(x') \in S_{x_1}-\{q\}$.
So $S_{x_1} = \{\phi(x'),q\}$ and $\phi(x') \not \in S_{x_2}-\{q\}$.

We extend $\phi$ by defining $\phi(x_1)=q$, defining $\phi(x_2)$ to be the element of $S_{x_2}-\{q\}$, and defining $\phi(x)$ to be an element of $S_x-\{q,\phi(x_2)\}$.
If $\phi(x)=\phi(y')$, then for every $z \in \{y_1,y_2\}$, defining $\phi(z)$ to be an element of $S_z-\{q'\}$ and defining $\phi(y)$ to be an element of $S_y-\{\phi(x),\phi(y_1),\phi(y_2)\}$ (which is possible since $\phi(x)=\phi(y') \not \in S_y$); if $\phi(x) \neq \phi(y')$, then defining $\phi(y_1)=\phi(y_2)=q'$ and defining $\phi(y)$ to be an element of $S_y-\{\phi(x),q'\}$.
Clearly, $\phi$ is a proper $L$-coloring.

Since $G$ is a counterexample, there exists a non-$R$-relaxed vertex $v$ in $G$ such that no color appears an odd number of times on $N_G(v)$.
Note that $v \in \{x,y\}$.
Since $\phi(x'),\phi(x_1)=q,\phi(x_2)$ are pairwise distinct, some color appears exactly once on $N_G(x)$.
So $v=y$.

By Claim 2 and the definition of $\phi$, $\phi(x)=\phi(y')$.
Hence $\phi(y_1) \neq \phi(y_2)$ by the definition of $\phi$.
Since $\phi(x)=\phi(y')$, we know that $\phi(y_1)$ or $\phi(y_2)$ appears exactly once on $N_G(y)$.
Hence $v \neq y$, a contradiction.
\end{pf}

\begin{lemma} \label{no_3_cycle_share}
No two distinct 3-cycles sharing edges. 
\end{lemma}

\begin{pf}
Since every 3-cycle cannot have $R$-length in $\{3,4,6\}$, it has exactly two edges in $R$, so it has $R$-length 5.
Since no two cycles of $R$-length 5 share exactly one edge, no two 3-cycles in $G$ share exactly one edge.
Since $G$ is simple, if two 3-cycles share edges, then they share exactly one edge.
This proves the lemma.
\end{pf}

\begin{lemma} \label{3_4_face_relaxed}
If $F_1$ is a 3-face and $F_2$ is a 4-face such that $F_1$ is adjacent to $F_2$, then $F_1$ and $F_2$ share exactly one edge, and all vertices incident with $F_1 \cup F_2$ are $R$-relaxed.
\end{lemma}

\begin{pf}
Since every vertex of $G$ has degree at least three (by Lemma \ref{min_deg}), and the length of $F_1$ and $F_2$ is at most four, we know that $F_1$ and $F_2$ are bounded by cycles $C_1$ and $C_2$, respectively.
Since $F_1$ is a 3-face and $F_2$ is a 4-face sharing at least one edge, they share exactly one edge by Lemma \ref{no_3_cycle_share}.
Since $C_1$ bounds a 3-face, it has $R$-length five by Lemma \ref{3_cycle_relaxed}.
Since $C_2$ shares exactly one edge with $C_1$, the $R$-length of $C_2$ is not in $\{3,4,5,6\}$.
So $|E(C_2)|+|R \cap E(C_2)| \geq 7$, implying that $|R \cap E(C_2)| \geq 3 = |E(C_2)|-1$.
Hence all vertices of $C_2$ are $R$-relaxed.
Recall that all vertices of $C_1$ are $R$-relaxed by Lemma \ref{3_cycle_relaxed}.
This proves the lemma.
\end{pf}

\begin{lemma}
Let $F_1$ be a 3-face, and let $F_2$ be a 4-face sharing exactly one edge with $F_1$.
Let $v$ be a 4-vertex incident with both $F_1$ and $F_2$.
Let $e$ be the edge incident with $F_1$ and $v$ but not incident with $F_2$.
Then the face incident with $e$ other than $F_1$ is a $(\geq 5)$-face. 
\end{lemma}

\begin{pf}
Let $F$ be the face incident with $e$ other than $F_1$.
By Lemma \ref{no_3_cycle_share}, $F$ has length at least four.
If $F$ has length four, then all vertices incident with $F_1 \cup F_2 \cup F$ are $R$-relaxed by Lemma \ref{3_4_face_relaxed}, so $v$ and all neighbors of $v$ are $R$-relaxed, contradicting Lemma \ref{relaxed_4_neighbors}.
\end{pf}

\begin{lemma} \label{3_not_on_4_4}
If $F_1$ and $F_2$ are 4-faces sharing exactly one edge, then every vertex incident with both $F_1$ and $F_2$ has degree at least four.
\end{lemma}

\begin{pf}
Since every vertex has degree at least three by Lemma \ref{min_deg}, $F_1$ and $F_2$ are bounded by 4-cycles $C_1$ and $C_2$, respectively.
Let $e$ be the edge shared by $C_1$ and $C_2$.
Note that the ends of $e$ are exactly the vertices incident with both $F_1$ and $F_2$ by Lemma \ref{no_3_cycle_share}.

If there exists $i \in [2]$ such that $C_i$ contains at least three edges in $R$, then all vertices in $C_i$ are $R$-relaxed, so every end of $e$ has at least two $R$-relaxed neighbors, and hence no end of $e$ has degree 3 by Lemma \ref{3_vertex_relaxed}.
So we may assume that each of $C_1$ and $C_2$ contains at most two edges in $R$.
Since $C_1$ and $C_2$ cannot have $R$-length 6, each of them has exactly one edge in $R$.
So there are cycles of $R$-length 5 sharing exactly one edge, a contradiction.
\end{pf}

\begin{lemma} \label{5_faces_share_3_neighbor}
Let $F_1$ and $F_2$ be 5-faces that share exactly one edge $a_1a_2$.
For each $i \in \{1,2\}$, let $a_i'$ be the neighbor of $a_i$ incident with $F_2$ other than $a_{3-i}$.
If $\deg(a_1)=3$, then $\deg(a_i') \geq 4$ for some $i \in \{1,2\}$.
\end{lemma}

\begin{pf}
Since every vertex has degree at least three by Lemma \ref{min_deg}, $F_1$ and $F_2$ are bounded by 5-cycles $C_1$ and $C_2$, respectively.

Suppose to the contrary that this lemma does not hold.
Then $\deg(a_1')=\deg(a_2')=3$ by Lemma \ref{min_deg}.
Since $a_1$ and $a_1'$ are adjacent 3-vertices, no vertex in $V(C_2)-\{a_1,a_1',a_2'\}$ can be $R$-relaxed by Lemma \ref{3_vertex_relaxed}.
So $E(C_2)-\{a_1a_1'\} \subseteq E(G)-R$.
Since $C_2$ cannot have $R$-length 6, $a_1a_1' \in E(G)-R$, so $C_2$ has $R$-length 5.

Similarly, since $a_1$ and $a_1'$ are adjacent 3-vertices, the two neighbors of $a_1$ in $C_1$ are not $R$-relaxed.
So $C_1$ contains at least four edges in $E(C)-R$.
Since $C_1$ cannot have $R$-length 6, $E(C_1)$ has $R$-length 5.
Then $C_1$ and $C_2$ are cycles of $R$-length 5 sharing exactly one edge, a contradiction.
\end{pf}

\begin{lemma} \label{no_3_v_554_f}
If $v$ is a 3-vertex incident with two 5-faces, then $v$ is not incident with a 4-face.
\end{lemma}

\begin{pf}
Suppose to the contrary that $v$ is incident with a 4-face $f$.
Let $f_1$ and $f_2$ be the 5-faces incident with $v$.
For each $i \in [2]$, let $e_i$ be the edge shared by $f$ and $f_i$ incident with $v$, and let $u_i$ be the end of $e_i$ other than $v$.
By Lemma \ref{3_vertex_relaxed}, $u_1$ or $u_2$ is non-$R$-relaxed.
By symmetry, we may assume that $u_1$ is non-$R$-relaxed.
So at most two edges incident with $f$ are in $R$.
Since the cycle bounding $f$ cannot have $R$-length 4 or 6, we know that exactly one edge incident with $f$ is in $R$, and the $R$-length of the cycle bounding $f$ is 5.
Hence $u_2$ is $R$-relaxed, and the cycle bounding $f_1$, denoted by $C_1$, does not have $R$-length 5.
By Lemma \ref{3_vertex_relaxed}, the neighbors of $v$ other than $u_2$ are non-$R$-relaxed.
So $C_1$ contains at least four edges not in $R$.
Since $C_1$ cannot have $R$-length 6, it has $R$-length 5, a contradiction.
\end{pf}

\begin{lemma} \label{no_4_v_4_f}
If $v$ is a 4-vertex, then some face incident with $v$ is not a 4-face.
\end{lemma}

\begin{pf}
Suppose to the contrary that $v$ is incident with four 4-faces.
Note that each 4-face is bounded by a 4-cycle by Lemma \ref{min_deg}.
Since each 4-cycle cannot have $R$-length 4 or 6, it contains one, three, or four edges in $R$.
So every 4-cycle either has $R$-length 5, or contains four $R$-relaxed vertices.

Since no two cycles of $R$-length 5 share exactly one edge, some 4-cycle bounding a 4-face incident with $v$ contains four $R$-relaxed vertices.
So $v$ is $R$-relaxed.

Let $v_1,v_2,v_3,v_4$ be the neighbors of $v$.
By Lemma \ref{relaxed_4_neighbors}, some $v_i$ is non-$R$-relaxed since $v$ is $R$-relaxed.
By symmetry, we may assume that $v_1$ is non-$R$-relaxed.
Let $C_1$ and $C_2$ be the cycles bounding the faces incident with $vv_1$.
Since $v_1$ is non-$R$-relaxed, each $C_i$ contains at least two edges not in $R$.
So the $R$-length of each $C_i$ equals 4, 5, or 6.
Since there is no cycle of $R$-length 4 or 6, both $C_1$ and $C_2$ have $R$-length 5.
But $C_1$ and $C_2$ share edges, so they share at least two edges.
However, it is impossible by Lemma \ref{min_deg}, a contradiction.
\end{pf}

\section{Discharging} \label{sec:discharging}

In this section, we assume that $G$ is a counterexample to Theorem \ref{main_real} such that $|V(G)|$ is as small as possible, and we fix a 2-cell embedding of $G$ in a surface of Euler genus at most $2$. 

For every $v \in V(G)$, let $\ch(v) = \deg(v)-4$.
For every face $f$, let $\ch(f) = \leng(f)-4$.

Now we apply the following discharging rules:

\begin{itemize}
    \item[(R1)] For every $(\geq 5)$-face $f$, $f$ sends $\frac{1}{2}$ unit of charge to every 3-vertex incident with $f$.
    \item[(R2)] For any $(\geq 5)$-face $f$, 3-face $f'$ and edge $ab$ shared by $f$ and $f'$ with $\deg(a)=\deg(b)=4$, if some edge incident with $f$ is between $a$ and a non-$R$-relaxed vertex in $N(a)-\{b\}$, and some vertex in $N(b)$ not incident with $f'$ is non-$R$-relaxed, then $f$ sends $\frac{1}{2}$ unit of charge to $f'$ via $ab$.
    \item[(R3)] For any $(\geq 5)$-face $f$, 3-face $f'$ and edge $ab$ shared by $f$ and $f'$ with $\deg(a)=\deg(b)=4$, if $f$ does not send charge to $f'$ via $ab$ due to (R2), then $f$ sends $\frac{1}{2}$ unit of charge to $f'$ via $ab$.
    \item[(R4)] For any $(\geq 5)$-face $f$ and edge $ab$ shared by $f$ and a 3-face $f'$ with $\deg(a)=4$ and $\deg(b) \geq 5$, $f$ sends $\frac{1}{4}$ unit of charge to $f'$ via $ab$.
    \item[(R5)] For any $(\geq 6)$-face $f$ and 5-face $f'$ sharing an edge $ab$ with $f$, if $a$ and $b$ are vertices such that $\deg(a)=3$, $\deg(b) =4$, and for each $v \in \{a,b\}$, there exists a unique vertex in $N(v)-\{a,b\}$ incident with $f'$ and this vertex is $R$-relaxed, then $f$ sends $\frac{1}{4}$ unit of charge to $f'$ via $ab$. 
    \item[(R6)] For any $(\geq 6)$-face $f$, any 4-vertex $v$ incident with $f$, and any edge $e$ incident with both $v$ and $f$, if 
        \begin{itemize}
            \item the end of $e$ other than $v$ is non-$R$-relaxed,
            \item some 5-face $f'$ shares $e$ with $f$, 
            \item some 5-face $f''$ distinct from $f$ shares exactly one edge $e'$ with $f'$ such that $v$ is incident with $e'$,
            \item the two neighbors of $v$ incident with $f''$ have degree three, and 
            \item $f''$ is adjacent to a 3-face, 
        \end{itemize}
        then $f$ sends $\frac{1}{4}$ unit of charge to $f''$ via $e'$.
    \item[(R7)] For every $(\geq 5)$-vertex $v$, $v$ sends $\frac{1}{2}$ unit of charge to each 3-face incident with $v$.
    \item[(R8)] For every $(\geq 6)$-vertex $v$, $v$ sends $\frac{1}{3}$ unit of charge to each 5-face $f$ incident with $v$ satisfying that none of the two edges incident with both $v$ and $f$ is incident with a 3-face.
\end{itemize}

\begin{figure}[h]
\centering
\renewcommand\thesubfigure{R\arabic{subfigure}}
\begin{subfigure}{0.20\textwidth}
     \centering
     \includegraphics[width=\textwidth]{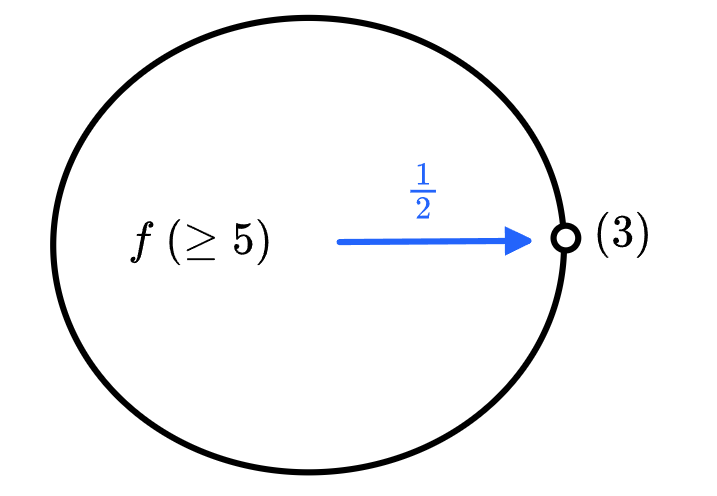}
     \caption{}
 \end{subfigure}
 ~
 \begin{subfigure}{0.24\textwidth}
     \centering
     \includegraphics[width=\textwidth]{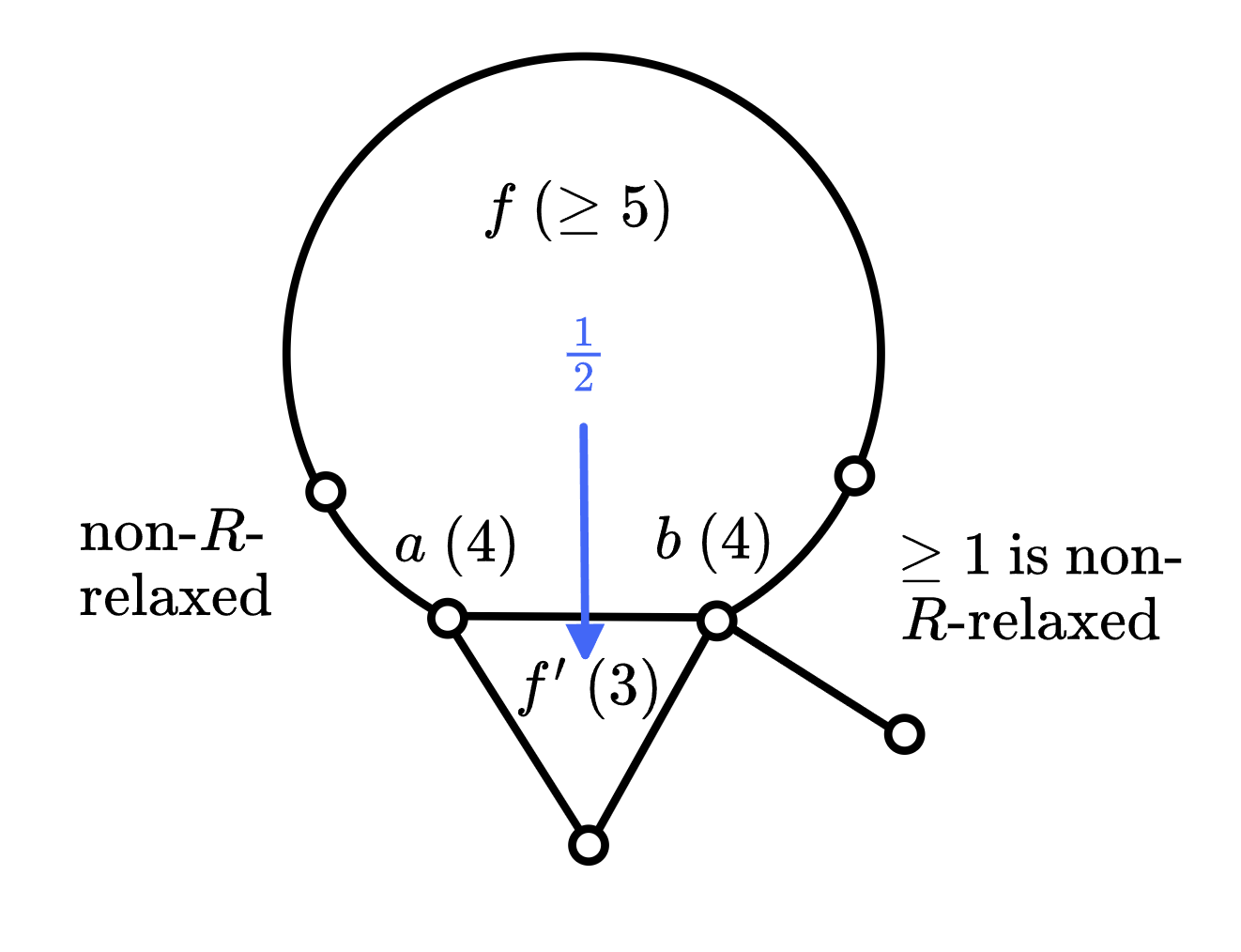}
     \caption{}
 \end{subfigure}
 ~
 \begin{subfigure}{0.22\textwidth}
     \centering
     \includegraphics[width=\textwidth]{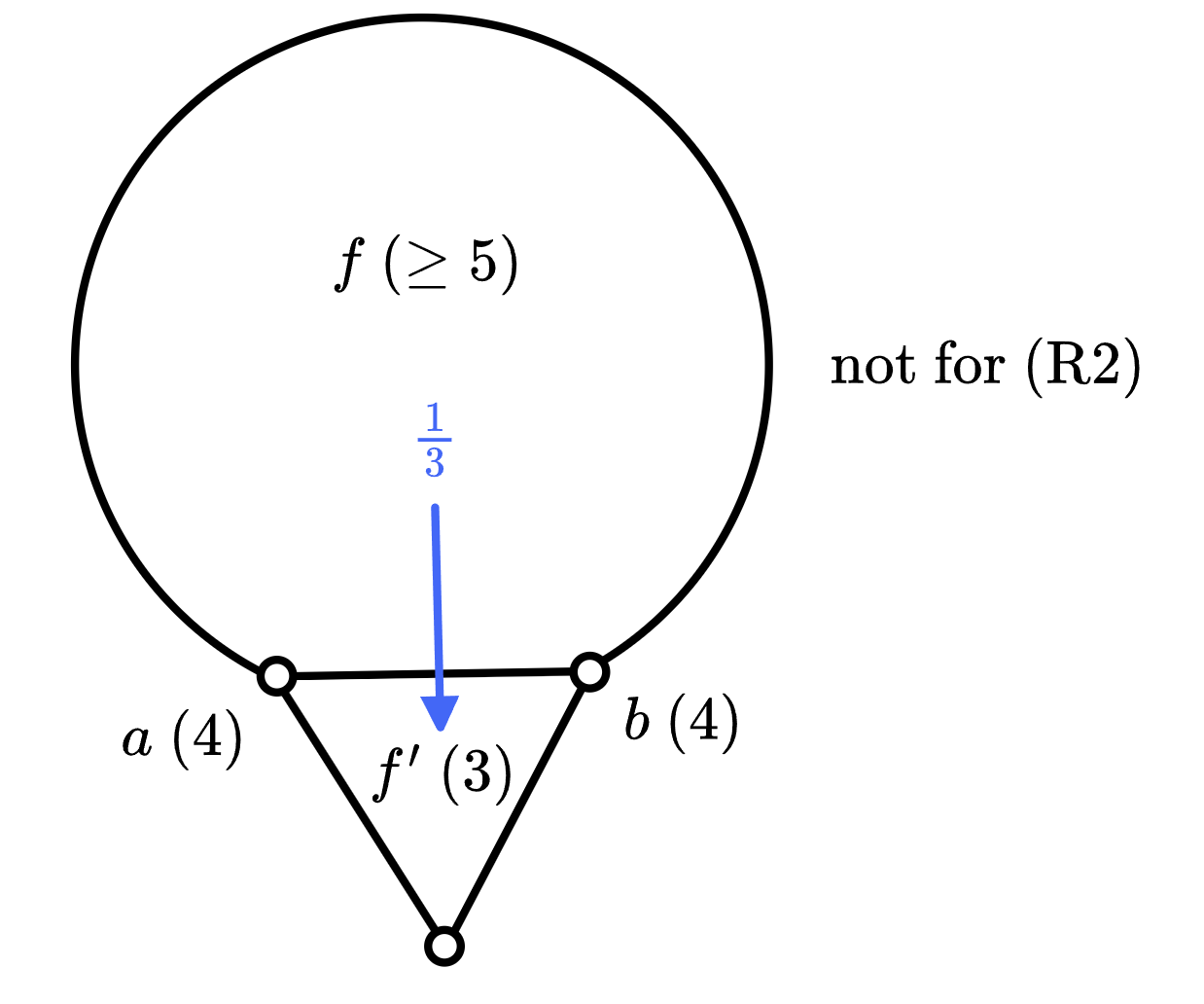}
     \caption{}
 \end{subfigure}
 ~
 \begin{subfigure}{0.2\textwidth}
     \centering
     \includegraphics[width=\textwidth]{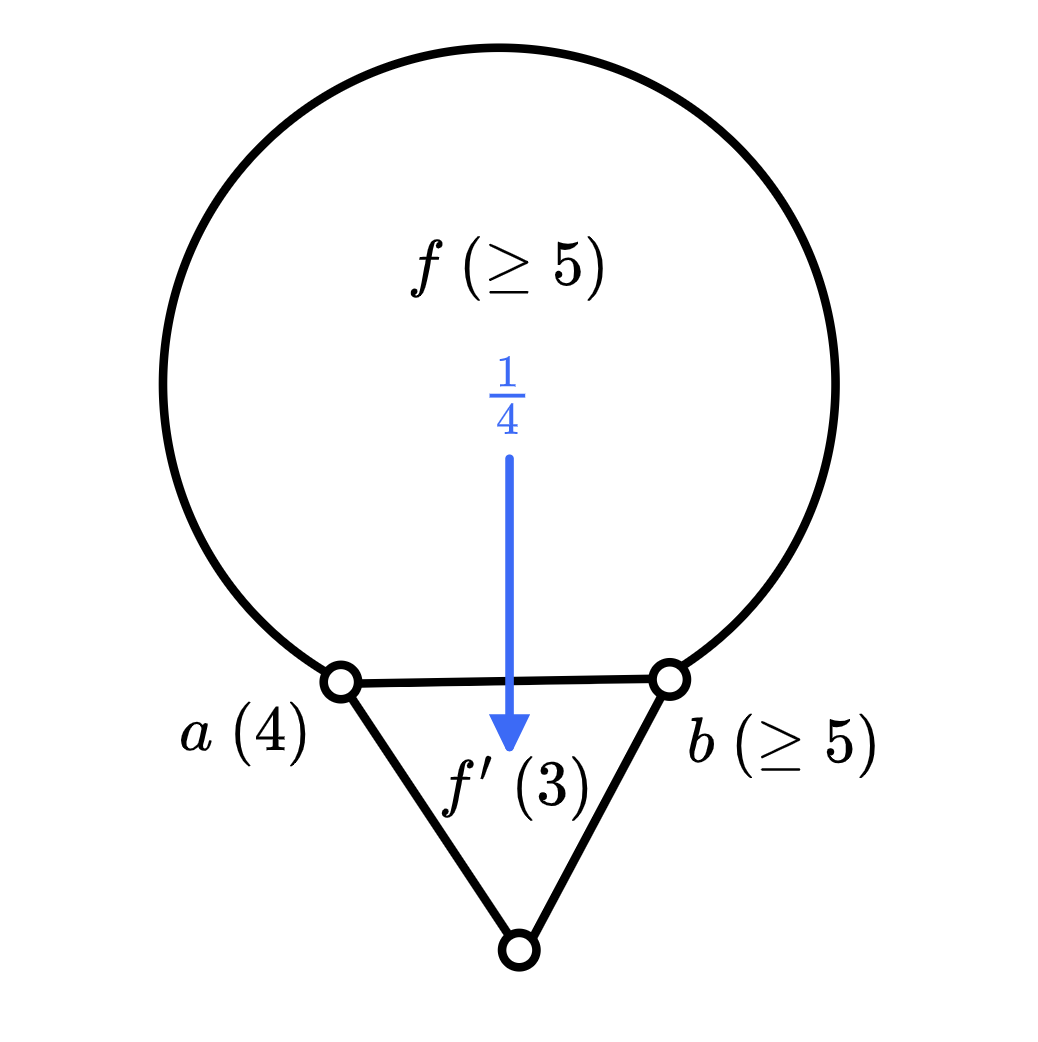}
     \caption{}
 \end{subfigure}
 \\
 \begin{subfigure}{0.24\textwidth}
     \centering
     \includegraphics[width=\textwidth]{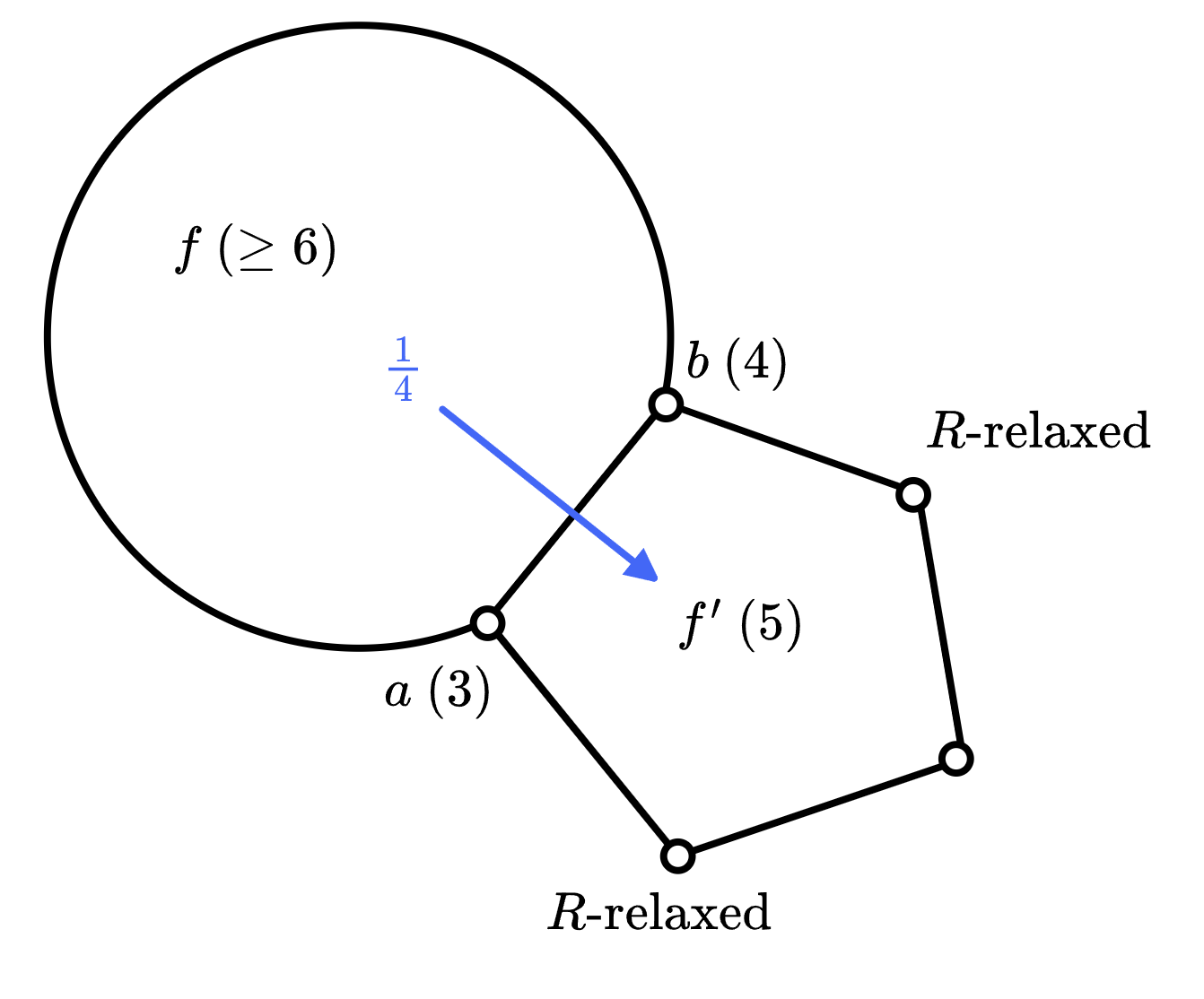}
     \caption{}
 \end{subfigure}
 ~
 \begin{subfigure}{0.22\textwidth}
     \centering
     \includegraphics[width=\textwidth]{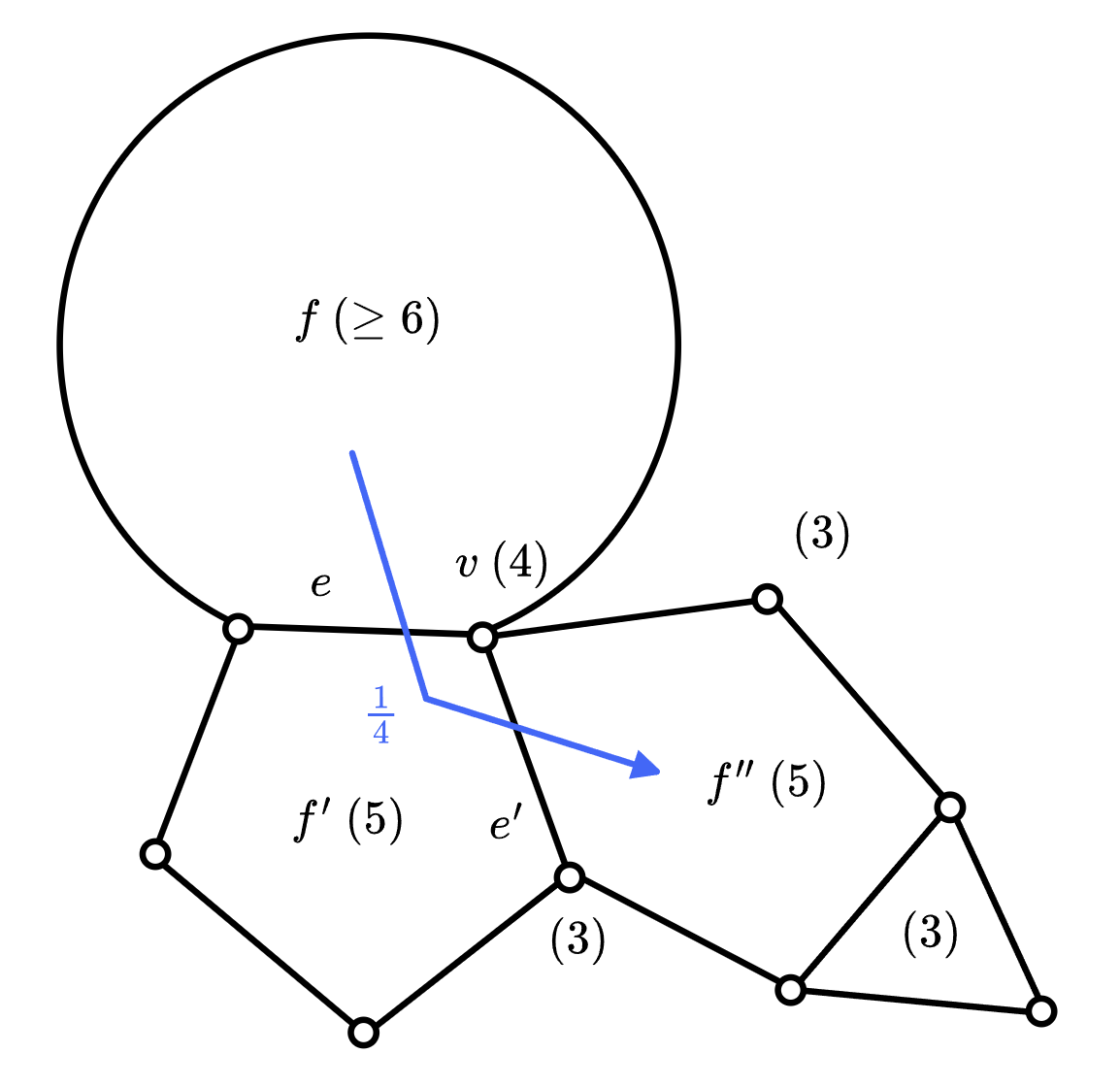}
     \caption{}
 \end{subfigure}
 ~
 \begin{subfigure}{0.21\textwidth}
     \centering
     \includegraphics[width=\textwidth]{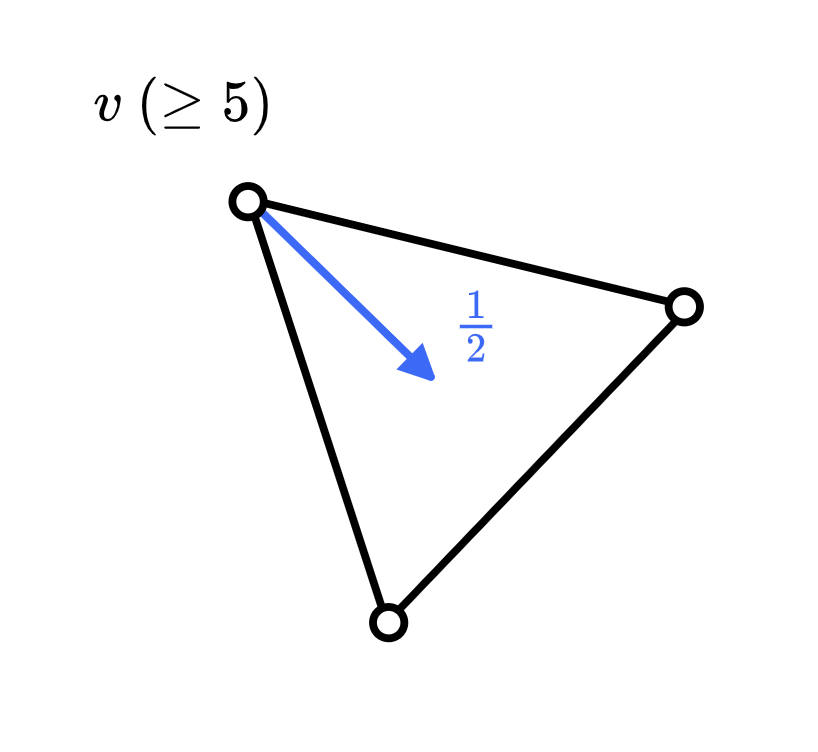}
     \caption{}
 \end{subfigure}
 ~
 \begin{subfigure}{0.24\textwidth}
     \centering
     \includegraphics[width=\textwidth]{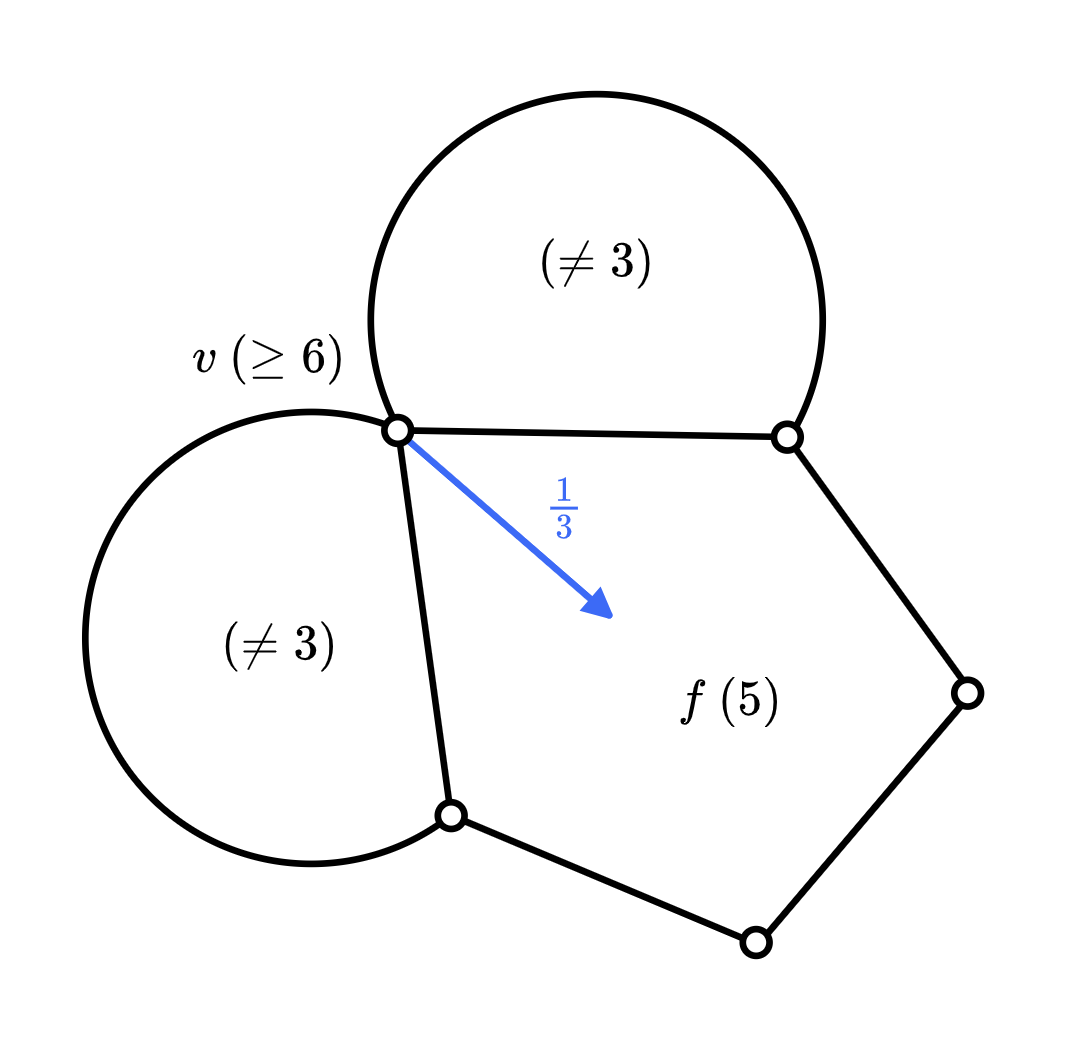}
     \caption{}
 \end{subfigure}

\caption{Diagrams for the rules (R1)-(R8). The numbers in the parentheses indicate the degree of a vertex or the length of a face. Note that two distinct faces possibly share more than one edge, so the actual situation could be more complicated than what the figures represent.}
\label{fig:discharging}
\end{figure}

For every vertex or face $x$ of $G$, let $\ch'(x)$ be the charge of $x$ after applying the above rules.

\begin{lemma} \label{discharging_vertex}
Let $v \in V(G)$.
Then $\ch'(v) \geq 0$. 
Moreover, if $\ch'(v)=0$, then either
    \begin{enumerate}
        \item $v$ is a 6-vertex and every face incident with $v$ is a 5-face, or
        \item $v$ is a 5-vertex and it is incident with exactly two 3-faces, or
        \item $v$ is a 4-vertex, or
        \item $v$ is a 3-vertex that is incident with exactly two $(\geq 5)$-faces.
    \end{enumerate}
\end{lemma}

\begin{pf}
Let $t$ be the number of 3-faces incident with $v$.
By Lemma \ref{no_3_cycle_share}, $t \leq \lfloor \frac{\deg(v)}{2} \rfloor$.

We first assume $\deg(v) \geq 6$.
Let $b$ be the number of 5-faces receiving charge from $v$ due to (R8).
Note that each 5-face receiving charge from $v$ due to (R8) does not share a 3-face that is incident with $v$.
So $b \leq (\deg(v)-t)-t = \deg(v)-2t$.
Hence by (R7) and (R8), $\ch'(v) = (\deg(v)-4) - \frac{t}{2} - \frac{b}{3} \geq \deg(v)-4 - \frac{t}{2} - \frac{\deg(v)-2t}{3} = \frac{2 \deg(v)}{3} -4 + \frac{t}{6} \geq \frac{t}{6} \geq 0$.
And if $\ch'(v)=0$, then $t=0$ and $b=\deg(v)=6$.

So we may assume $\deg(v) \leq 5$.

If $v$ is a 5-vertex, then $t \leq \lfloor \frac{\deg(v)}{2} \rfloor = 2$, so by (R7), $\ch'(v) = (\deg(v)-4) - \frac{t}{2} \geq 0$; and the equality holds when $t=2$.

If $v$ is a 4-vertex, then $\ch'(v)=\ch(v)=\deg(v)-4=0$.

Let $k$ be the number of $(\geq 5)$-faces incident with $v$.
If $v$ is a 3-vertex, then $k \geq 2$ by Lemmas \ref{3_not_3} and \ref{3_not_on_4_4}, so $\ch'(v) = (\deg(v)-4)+ \frac{k}{2} \geq 0$ by (R1); the equality holds when $k=2$.
\end{pf}

\begin{lemma} \label{discharging_3_face}
Let $f$ be a 3-face.
Then $\ch'(f) \geq 0$.
Moreover, if $\ch'(f)=0$, then either 
    \begin{enumerate}
        \item $f$ is incident with a 5-vertex and two 4-vertices such that either
            \begin{enumerate}
                \item the face sharing the edge incident with $f$ formed by the 4-vertices is a 4-face and the faces sharing the other two edges incident with $f$ are $(\geq 5)$-faces, or
                \item the face sharing the edge incident with $f$ formed by the 4-vertices is a $(\geq 5)$-face and the faces sharing the other two edges incident with $f$ are $4$-faces, or
            \end{enumerate}
        \item $f$ is incident with three 4-vertices such that exactly one 4-face is adjacent to $f$. 
    \end{enumerate}
\end{lemma}

\begin{pf}
By Lemma \ref{no_3_cycle_share}, all faces adjacent to $f$ are $(\geq 4)$-faces.
By Lemma \ref{3_not_3}, all vertices incident with $f$ are $(\geq 4)$-vertices.

If all vertices incident with $f$ are $(\geq 5)$-vertices, then $\ch'(f) \geq (3-4)+3 \cdot \frac{1}{2} >0$ by (R7).
So we may assume that $f$ is incident with at most two $(\geq 5)$-vertices.
Hence $f$ is incident with a 4-vertex $a$.

Denote the vertices incident with $f$ other than $a$ by $b$ and $c$.

Let $f_b$ and $f_c$ be the faces sharing $ab$ and $ac$ with $f$, respectively.
Let $a'$ be a vertex in $N(a)-\{b,c\}$ incident with $f_b$, and let $a''$ be a vertex in $N(a)-\{b,c\}$ incident with $f_c$ such that $N(a)=\{b,c,a',a''\}$.
By Lemma \ref{3_cycle_relaxed}, all $a,b,c$ are $R$-relaxed.
By Lemma \ref{relaxed_4_neighbors}, $a'$ or $a''$ is non-$R$-relaxed.
By the symmetry between $b$ and $c$, we may assume that $a'$ is non-$R$-relaxed.
So $f_b$ is a $(\geq 5)$-face by Lemma \ref{3_4_face_relaxed}.

Hence if $b$ and $c$ are $(\geq 5)$-vertices, then $f_b$ sends $\frac{1}{4}$ unit of charge to $f$ due to (R4), and each $b$ and $c$ sends $\frac{1}{2}$ unit of charge to $f$ due to (R7), so $\ch'(f) \geq (3-4) + \frac{1}{4} + 2 \cdot \frac{1}{2} >0$ and the lemma holds.

So we may assume that at least one of $b$ and $c$ is a 4-vertex.
Let $f'$ be the face sharing $bc$ with $f$.

We first assume $\deg(b) \geq 5$.
So $c$ is a 4-vertex.
By Lemmas \ref{3_cycle_relaxed}, \ref{relaxed_4_neighbors} and \ref{3_4_face_relaxed}, at least one of $f'$ and $f_c$ is a $(\geq 5)$-face.
Note that $f_b$ sends $\frac{1}{4}$ unit of charge to $f$ via $ab$ due to (R4), and $b$ sends $\frac{1}{2}$ unit of charge to $f$ due to (R7).
If $f_c$ is a $(\geq 5)$-face, then $f_c$ sends at least $\frac{1}{3}$ unit of charge to $f$ via $ac$ due to (R2) or (R3), so $\ch'(f) \geq (3-4) + \frac{1}{4} + \frac{1}{2} + \frac{1}{3} > 0$.
If $f_c$ is a 4-face, then $f'$ is a $(\geq 5)$-face sending $\frac{1}{4}$ unit of charge to $f$ via $bc$ due to (R4), so $\ch'(f) \geq (3-4) + \frac{1}{4} + \frac{1}{2} + \frac{1}{4}=0$; when $\ch'(f)=0$, Statement 1(a) of this lemma holds.

So we may assume $\deg(b)=4$.
By Lemmas \ref{3_cycle_relaxed} and \ref{relaxed_4_neighbors}, some vertex in $N(b)$ not incident with $f$ is non-$R$-relaxed.
Recall that $a'$ is non-$R$-relaxed, so $f_b$ sends $\frac{1}{2}$ to $f$ via $ab$ due to (R2).
If $\deg(c) \geq 5$, then $c$ sends $\frac{1}{2}$ unit of charge to $f$ due to (R7), so $\ch'(f) \geq (3-4)+ \frac{1}{2}+\frac{1}{2}=0$; the equality holds only when $f'$ and $f_c$ do not send charge to $f$ via $bc$ and $ac$, respectively, implying both $f'$ and $f_c$ are 4-faces by (R4), so Statement 1(b) of this lemma holds.

So we may assume $\deg(c)=4$.
By Lemmas \ref{3_cycle_relaxed} and \ref{relaxed_4_neighbors}, some neighbor $c'$ of $c$ not incident with $f$ is non-$R$-relaxed.
By Lemma \ref{3_4_face_relaxed}, the face $f''$ adjacent to $f$ incident with $cc'$ is a $(\geq 5)$-face.
Let $x \in \{a,b,c\}$ such that $f$ and $f''$ share $cx$.
By (R2), $f''$ sends $\frac{1}{2}$ unit of charge to $f$ via $xc$.
So $\ch'(f) \geq (3-4) + \frac{1}{2} + \frac{1}{2} =0$; the equality holds only when the face sharing the unique edge $e \in \{ab,ac,bc\}-\{ab,xc\}$ does not send charge to $f$ via $e$, implying that this face is a 4-face by (R2) and (R3), so Statement 2 of this lemma holds.
This proves the lemma.
\end{pf}

\begin{lemma} \label{discharging_4_face}
$\ch'(f) \geq 0$ for every $4$-face $f$.
\end{lemma}

\begin{pf}
Note that $\ch(f)=0$.
So this lemma is obvious since every $4$-face does not receive or send charge.
\end{pf}

\begin{lemma} \label{discharging_6_face}
Let $f$ be a $(\geq 6)$-face.
Then $\ch'(f) \geq 0$.
Moreover, if $\ch'(f)=0$, then $f$ is a 6-face bounded by a 6-cycle such that either
    \begin{enumerate}
        \item $f$ is incident with at least one $(\geq 6)$-vertex, or
        \item the degree of the vertices incident with $f$ are 3, 3, 4, 3, 3, 4 listed in the cyclic order given by the boundary cycle of $f$.
    \end{enumerate}
\end{lemma}

\begin{pf}
Let $W$ be a walk that defines the length of $f$.
By Lemma \ref{conn}, $G$ is connected, so $W$ is connected.
Note that $W$ possibly contains repeated edges.

During the proof, we will send charge from $f$ required by rules (R1)-(R6) in several steps; in the meantime, we will put blue and green tokens on edges in $W$ to indicate how much charge has been sent from $f$ so far.
We denote the number of blue tokens and green tokens by $b$ and $g$, respectively, and show that at any moment of the process, the following two conditions are maintained:
    \begin{itemize}
        \item[(i)] The charge that has been sent from $f$ is at most $\frac{1}{4}b+\frac{1}{3}g$.
        \item[(ii)] If an edge receives at least two tokens, then it receives exactly two blue tokens and no green tokens.
    \end{itemize}

Whenever we send charge via an edge $e$ due to (R3) or (R4), we put a green token on $e$.
Clearly, (i) and (ii) hold at this point.

Whenever we send charge via an edge $e$ due to (R5), we put a blue token on $e$.
Clearly, (i) and (ii) hold at this point.
In fact, no edge receives at least two tokens at this point.

Whenever we send charge via an edge $e'$ due to (R6), we put a blue token on $e$, where $e$ is the edge $e$ stated in (R6).
Clearly, (i) holds and no edge receives both blue and green tokens at this point.
Moreover, if an edge receives at least two blue tokens, then this edge receives exactly two tokens, and both tokens come from applying (R6) (since no edge receiving tokens due to (R6) also received tokens due to (R3)-(R5)); but since every vertex in a 3-cycle is $R$-relaxed, the ends of this edge are 4-vertices each adjacent to two 3-vertices that have relaxed neighbors, contradicting Lemma \ref{relaxed_443333}.
So no edge receives at least two tokens and hence (ii) holds.

For every edge $ab$ that $f$ sends charge via $ab$ due to (R2), we know that there exists an edge $e_a$ incident with $f$ between $a$ and a non-$R$-relaxed vertex $a'$ in $N(a)-\{b\}$; note that $e_a$ does not receive any token due to (R3)-(R5) since $a'$ is non-$R$-relaxed and every vertex in a 3-cycle is $R$-relaxed by Lemma \ref{3_cycle_relaxed}, and $e_a$ does not receive any token due to (R6) since $a$ is $R$-relaxed (by Lemma \ref{3_cycle_relaxed}) and $a$ has two neighbors in 3-cycles and those two neighbors cannot be 3-vertices by Lemma \ref{3_not_3}.
Whenever $f$ sends charge via $ab$ due to (R2), we put a blue token on each $ab$ and $e_a$.
Clearly, (i) holds at this point.
Moreover, this step does not create any edge that has two tokens.
So no edge received at least two tokens and (ii) holds.

Note that for every edge that has received a token, at most one of its ends is a 3-vertex.
Then we send $\frac{1}{2}$ charge to each 3-vertex $v$ due to (R1), select two edges incident with both $v$ and $f$ such that for each selected edge,
    \begin{itemize}
        \item its both ends are 3-vertices if possible, and 
        \item subject to this, $f$ is the unique face incident with it if possible,
    \end{itemize} 
and put a blue token on each of those selected edges. 
Clearly, (i) holds and no edge received both blue tokens and green tokens.

For every edge $e$ incident with $f$, let $b_e$ and $g_e$ be the number of tokens on $e$, respectively, and let $w_e = \frac{1}{4} \cdot b_e + \frac{1}{3} \cdot g_e$.
Note that (i) implies that the charge sent from $f$ equals $\sum_e w_e$, where the sum is over all edges $e$ incident with $f$ (with multiplicity).

Note that if $e$ is an edge that has at least two blue tokens, then either both ends of $e$ are 3-vertices, or $e$ receives a blue token due to (R5) and a blue token due to (R1).
So (ii) holds.
Moreover, if both ends of some edge $e$ are 3-vertices, then by Lemma \ref{3_vertex_relaxed} and the rule for selecting edges to put blue tokens due to (R1), we know that no edge incident with $e$ received two tokens and no edge incident with $e$ received green tokens, and we let $\iota(e)$ be an edge different from $e$ but incident with an end of $e$. 
If some edge $ab$ with $\deg(a)=3$ and $\deg(b)=4$ receives a blue token due to (R5), then by Lemmas \ref{3_vertex_relaxed} and \ref{relaxed_433}, some edge $q_{ab} \neq ab$ incident with both $f$ and $b$ does not receive two blue tokens and does not receive green tokens; and we can choose $\iota(ab)$ such that there exists at most one edge $a'b'$ different from $ab$ receiving two blue tokens such that $\iota(ab)=\iota(a'b')$. 
Hence $\iota$ is a function such that the preimage of any element of the codomain of $\iota$ has size at most two, and every element of the image of $\iota$ does not receive a two blue tokens and does not receive green tokens.
So at most $\lfloor \frac{2}{3} \cdot \leng(f) \rfloor$ edges incident with $f$ receive at least two tokens, and if exactly $\lfloor \frac{2}{3} \cdot \leng(f) \rfloor$ edges incident with $f$ receive at least two tokens, then some edge does not receive a green token and does not receive two blue tokens. 
That is, 
    \begin{itemize}
        \item $w_e = \frac{1}{4} \cdot b_e + \frac{1}{3} \cdot g_e \leq \frac{1}{2}$ for every edge $e$ incident with $f$,
        \item all but at most $\lfloor \frac{2}{3} \cdot \leng(f) \rfloor$ edges $e$ incident with $f$ satisfy $w_e \leq \frac{1}{3}$,
        \item if there are exactly $\lfloor \frac{2}{3} \cdot \leng(f) \rfloor$ edges $e$ satisfy $w_e > \frac{1}{3}$, then $w_{e'} < \frac{1}{3}$ for some edge $e'$ incident with $f$.
    \end{itemize}
Therefore, the charge sent from $f$ is strictly less than $\frac{1}{3} \cdot \leng(f) + (\frac{1}{2} - \frac{1}{3}) \cdot \lfloor \frac{2}{3} \cdot \leng(f) \rfloor \leq \frac{4}{9} \cdot \leng(f)$; when $\leng(f)=7$, $f$ sends strictly less than $\frac{1}{3} \cdot \leng(f) + (\frac{1}{2} - \frac{1}{3}) \cdot \lfloor \frac{2}{3} \cdot \leng(f) \rfloor = \frac{1}{3} \cdot 7 + \frac{1}{6} \cdot 4 = 3$.
So if $\leng(f) \geq 8$, then $\ch'(f) > (\leng(f)-4)-\frac{4}{9} \cdot \leng(f) > 0$; if $\leng(f)=7$, then $\ch'(f) > (\leng(f)-4)-3 = 0$.

Hence, this lemma holds when $\leng(f) \geq 7$.
So we may assume $\leng(f)=6$.

If $W$ is not a cycle, then since $G$ has minimum degree at least three, $W$ is a union of two 3-cycles that share exactly one vertex, so none of (R1)-(R6) can apply by Lemmas \ref{3_not_3}, \ref{no_3_cycle_share} and \ref{relaxed_4_neighbors}, and hence $\ch'(f)=\ch(f)>0$.
So we may assume that $W$ is a 6-cycle.

Denote $W$ by $v_1v_2...v_6v_1$.

Suppose that $f$ does not satisfy the conclusion of this lemma.
In particular, $\ch'(f) \leq 0$ and, if $\ch'(f)=0$, then there is no $(\geq 6)$-vertex incident with $f$.
We shall derive a contradiction to complete the proof.

\medskip

\noindent{\bf Claim 1:} There does not exist $i \in [6]$ such that $v_iv_{i+1},v_{i+1}v_{i+2},v_{i+2}v_{i+3}$ have green tokens, where the indices are computed modulo 6.

\noindent{\bf Proof of Claim 1:}
Suppose to the contrary.
Then $v_{i+1}$ or $v_{i+2}$ is a 4-vertex such that it and all its neighbors are contained in 3-cycles, contradicting Lemmas \ref{3_cycle_relaxed} and \ref{relaxed_4_neighbors}.
$\Box$

\medskip

\noindent{\bf Claim 2:} $f$ does not send charge due to (R5).

\noindent{\bf Proof of Claim 2:}
Suppose to the contrary that $f$ sends charge due to (R5).
By symmetry, we may assume that $f$ sends charge via $v_2v_3$ due to (R5), and assume that $\deg(v_2)=3$ and $\deg(v_3)=4$.
Since $W$ is a cycle, each $v_3$ and $v_4$ is adjacent to an $R$-relaxed vertex not incident with $f$.
By Lemma \ref{3_vertex_relaxed}, $v_1$ is a $(\geq 4)$-vertex and the edges $v_6v_1,v_1v_2,v_2v_3,v_3v_4$ are not in $R$.
The same argument shows that if there exists an edge $e$ in $W$ other than $v_1v_2$ and $v_2v_3$ such that $f$ sends charge via $e$ due to (R5), then since $v_1$ is a $(\geq 4)$-vertex, $v_4v_5$ and $v_5v_6$ are not in $R$, so all edges of $W$ are not in $R$, and hence $W$ is a cycle of $R$-length 6, a contradiction.

So $v_1v_2$ is the only possible edge other than $v_2v_3$ for which $f$ sends charge via due to (R5).

We first suppose that $f$ sends charge via $v_1v_2$ due to (R5).
Then $v_1$ is a 4-vertex, and by Lemma \ref{relaxed_433}, no two 3-vertices in $W$ are adjacent.
So $v_1v_2$ and $v_2v_3$ are the only edges received two blue tokens.
That is, $w_{v_1v_2}=w_{v_2v_3}=\frac{1}{2}$ and $w_e \leq \frac{1}{3}$ for every $e \in E(W)-\{v_1v_2,v_2v_3\}$.
Moreover, $v_6v_1$ and $v_3v_4$ do not receive green tokens by Lemma \ref{relaxed_4_neighbors}, so $w_e \leq \frac{1}{4}$ for every $e \in \{v_1v_6,v_3v_4\}$.
Hence, if $v_3v_4$ does not have a blue token, then $v_3v_4$ does not have a token and $\sum_{e \in E(W)}w_e \leq 2 \cdot \frac{1}{2} + (0+\frac{1}{4}) + 2 \cdot \frac{1}{3} < 2$, so $\ch'(f)>0$, a contradiction.
So $v_3v_4$ has a blue token.
Similarly, $v_1v_6$ has a blue token.
By Lemma \ref{3_vertex_relaxed}, $v_1$ and $v_3$ are non-$R$-relaxed, so $v_6v_1,v_1v_2,v_2v_3,v_3v_4$ are not in $R$.
Since $W$ cannot have $R$-length 6, $v_4v_5$ or $v_5v_6$ is in $R$.
It implies that $v_5$ is $R$-relaxed, so $v_4$ or $v_6$ are not 3-vertices by Lemma \ref{relaxed_433}.
By symmetry, we may assume that $v_4v_5 \in R$, so $v_4$ and $v_5$ are $R$-relaxed.
Hence $v_3v_4$ does not receive a blue token due to (R6) by Lemma \ref{relaxed_433}.
Since $v_3v_4$ has a blue token, $f$ sends charge via $v_4v_5$ due to (R2).
So $\deg(v_5)=4$ and some neighbor of $v_5$ not incident with the face sharing $v_4v_5$ with $f$ is non-$R$-relaxed.
Hence $v_5v_6$ cannot have a green token by Lemmas \ref{3_cycle_relaxed} and \ref{relaxed_4_neighbors}.
Since $v_5$ is $R$-relaxed, $v_5v_6$ does not receive a blue token due to (R6) by Lemma \ref{3_not_3}.
Since $v_5$ is $R$-relaxed and $v_5$ and $v_6$ are not 3-vertices, $f$ cannot send charge via $v_5v_6$ or $v_6v_1$ due to (R2) or (R5).
So $v_5v_6$ has no token.
Hence $\sum_{e \in E(W)}w_e \leq 2 \cdot \frac{1}{2} + 2 \cdot \frac{1}{4} + 0 + \frac{1}{4} < 2$ and $\ch'(f)>0$, a contradiction.

So $f$ does not send charge via $v_1v_2$ due to (R5).
Hence $f$ does not send charge via an edge other than $v_2v_3$ due to (R5).
Since $v_2$ is a 3-vertex and $v_1v_2$ does not have two blue tokens, $w_{v_1v_2}=\frac{1}{4}$.

Now we suppose that some edge in $W$ other than $v_2v_3$ have two blue tokens.
Then both ends of this edge are 3-vertices.
Since $f$ sends charge via $v_2v_3$ due to (R5), this edge must be $v_5v_6$ by Lemmas \ref{3_vertex_relaxed} and \ref{relaxed_433}.
Then $v_4v_5$ and $v_6v_1$ do not have green tokens and $v_3v_4$ does not have any token according to Lemma \ref{relaxed_433}.
Hence $\sum_{e \in E(W)}w_e \leq 2 \cdot \frac{1}{2} + 3 \cdot \frac{1}{4} < 2$, so $\ch'(f)>0$, a contradiction.

So $v_2v_3$ is the only edge that has two blue tokens.
By Lemma \ref{relaxed_4_neighbors}, $v_3v_4$ does not have a green token, so $w_{v_3v_4} \leq \frac{1}{4}$.
Hence $w_{v_1v_2}+w_{v_2v_3}+w_{v_3v_4} \leq \frac{1}{2} + 2 \cdot \frac{1}{4}=1$.
By Claim 1, at least one of $v_4v_5,v_5v_6,v_6v_1$ does not receive a green token, and if two of them have green tokens, then the remaining one does not have a blue token by Lemma \ref{relaxed_433}, so $w_{v_4v_5}+w_{v_5v_6}+w_{v_6v_1} \leq \max\{2 \cdot \frac{1}{3}, \frac{1}{3} + 2 \cdot \frac{1}{4}\} <1$.
Hence $\sum_{e \in E(W)}w_e < 2$, so $\ch'(f)>0$, a contradiction.
$\Box$

\medskip

\noindent{\bf Claim 3:} No edge received two blue tokens.

\noindent{\bf Proof of Claim 3:}
Suppose to the contrary that some edge, say $v_2v_3$, receives two blue tokens.
By Claim 2, both $v_2$ and $v_3$ are 3-vertices.
So $v_1$ and $v_4$ are non-$R$-relaxed vertices by Lemma \ref{3_vertex_relaxed}.
In particular, $v_1$ and $v_4$ have even degree. 

We first suppose that there exists an edge other than $v_2v_3$ having two blue tokens.
Then this edge must be $v_5v_6$, and $v_5$ and $v_6$ are 3-vertices by Claim 2.
For every edge $e \in E(W)-\{v_2v_3,v_5v_6\}$, since $e$ has a blue token, (ii) implies that $w_e=\frac{1}{4}$.
So $\ch'(f) = (\leng(f)-4)-\sum_{e \in E(W)}w_e = 2-(4 \cdot \frac{1}{4} + 2 \cdot \frac{1}{2}) = 0$.
Hence $W$ does not contain a $(\geq 6)$-vertex.
Since $v_1$ and $v_4$ have even degree, we know $\deg(v_1)=\deg(v_4)=4$.
Hence Statement 2 of this lemma holds, a contradiction.

So no edge other than $v_2v_3$ has two blue tokens.
Since both $v_2$ and $v_3$ are 3-vertices, $v_1v_2$ and $v_3v_4$ have no green tokens, so $w_{v_1v_2}+w_{v_3v_4}=\frac{1}{2}$.
By Claim 1, at least one of $v_4v_5,v_5v_6,v_6v_1$ does not have a green token, so $w_{v_4v_5}+w_{v_5v_6}+w_{v_6v_1} <1$.
Hence $\sum_{e \in E(W)}w_e < \frac{1}{2}+\frac{1}{2}+1 = 2$, so $\ch'(f)>0$, a contradiction.
$\Box$

\medskip

By Claim 3, $w_e \leq \frac{1}{3}$ for every $e \in E(W)$.
By Claim 1, $w_e < \frac{1}{3}$ for some $e \in E(W)$.
So $\sum_{e \in E(W)}w_e < 6 \cdot \frac{1}{3} = 2$.
Hence $\ch'(f)>0$, a contradiction.
This proves the lemma.
\end{pf}

\begin{lemma} \label{discharging_5_face}
Let $f$ be a 5-face.
Then $\ch'(f) \geq 0$.
Moreover, if $\ch'(f)=0$, then $f$ is bounded by a 5-cycle $v_1v_2v_3v_4v_5v_1$ such that either
    \begin{enumerate}
        \item $f$ is incident with at least two 3-vertices, or
        \item $f$ is incident with at least one 3-face and at least one 3-vertex, or
        \item there exists $i \in [5]$ such that 
            \begin{enumerate}
                \item $v_iv_{i+1}$ is incident with a 3-face $f_i$, 
                \item $v_{i+2}v_{i+3}$ is incident with a 3-face $f_{i+2}$, 
                \item both $f_i$ and $f_{i+2}$ receive charge from $f$ due to (R2), 
                \item the neighbor of $v_{i+1}$ not incident with $f$ or $f_i$ is non-$R$-relaxed, and
                \item the neighbor of $v_{i+2}$ not incident with $f$ or $f_{i+2}$ is non-$R$-relaxed,
                \item $v_{i+4}$ is non-$R$-relaxed,
            \end{enumerate}
            where the indices are computed modulo $5$.
    \end{enumerate}
\end{lemma}

\begin{pf}
Since $f$ is a 5-face and $G$ has minimum degree at least three, $f$ is bounded by a 5-cycle $v_1v_2v_3v_4v_5v_1$, denoted by $C$.
For every $i \in [5]$, let $f_i$ be the face sharing $v_iv_{i+1}$ with $f$, where $v_6=v_1$.
Since $G$ has minimum degree at least three, if $f_i$ is a 3-face, then the vertex incident with $f_i$ other than $v_i$ and $v_{i+1}$ is not incident with $f$.

Suppose that $f$ is a counterexample to this lemma.
That is, either $\ch'(f)<0$, or $\ch'(f)=0$ and none of Statements 1-3 holds.

\medskip

\noindent{\bf Claim 1:} If there are two 3-faces receiving charge from $f$ due to (R2), (R3), or (R4), then the triangles bounding these two faces do not share a vertex incident with $f$.

\noindent{\bf Proof of Claim 1:}
Suppose to the contrary.
Without loss of generality, we may assume that $f_1$ and $f_2$ are 3-faces receiving charge from $f$ due to (R2), (R3), or (R4).
By Lemmas \ref{3_cycle_relaxed} and \ref{no_3_cycle_share}, $v_2$ is an $R$-relaxed vertex that has four $R$-relaxed neighbors.
So $v_2$ is a $(\geq 5)$-vertex by Lemma \ref{relaxed_4_neighbors}.
Hence $f_1$ and $f_2$ do receive charge from $f$ due to (R2) or (R3).
Since $f_1$ and $f_2$ receive charge from $f$ due to (R2), (R3), or (R4), we know that they receive charge due to (R4), and $v_1$ and $v_3$ are $4$-vertices.
So $f$ sends $\frac{1}{4}$ to each of $f_1$ and $f_2$.
Since $v_1$ and $v_3$ are 4-vertices contained in triangles, Lemmas \ref{3_cycle_relaxed} and \ref{relaxed_4_neighbors} imply that $f_5$ and $f_3$ are not 3-faces.
So if $f_4$ is not a 3-face and $v_4$ and $v_5$ are $(\geq 4)$-vertices, then $\ch'(f) \geq (5-4) - 2 \cdot \frac{1}{4}>0$, a contradiction.

Hence either $f_4$ is a 3-face, or at least one of $v_4$ and $v_5$ is a 3-vertex.

We first suppose that $f_4$ is a 3-face.
By Lemma \ref{3_not_3}, $v_4$ and $v_5$ are $(\geq 4)$-vertices.
So no 3-vertex is incident with $f$.
Since $v_1$ and $v_3$ are contained in triangles, they are $R$-relaxed by Lemma \ref{3_cycle_relaxed}.
Hence $f$ does not send charge to $f_4$ due to (R2).
So if $f$ sends charge to $f_4$, then it sends at most $\frac{1}{3}$ to $f_4$.
Therefore, $\ch'(f) \geq (5-4) - 2 \cdot \frac{1}{4} - \frac{1}{3} >0$, a contradiction.

So $f_4$ is not a 3-face.
Hence $v_4$ or $v_5$ is a 3-vertex.
Since $v_1$ is $R$-relaxed, it is impossible that both $v_4$ and $v_5$ are 3-vertices by Lemma \ref{3_vertex_relaxed}.
So exactly one of $v_4$ and $v_5$ is a 3-vertex.
Hence $\ch'(f) = (5-4) - 2 \cdot \frac{1}{4} - \frac{1}{2}=0$, so Statement 2 of this lemma holds, contradicting that $f$ is a counterexample.
$\Box$

\medskip

\noindent{\bf Claim 2:} There do not exist two 3-faces receiving charge from $f$ due to (R2), (R3), or (R4).

\noindent{\bf Proof of Claim 2:}
Suppose to the contrary that there exist two 3-faces receiving charge from $f$ due to (R2), (R3), or (R4).
By Claim 1, they do not share a vertex incident with $f$.
So we may without loss of generality assume that $f_1$ and $f_3$ are 3-faces receiving charge from $f$ due to (R2), (R3), or (R4).
By Lemma \ref{3_cycle_relaxed}, $v_1,v_2,v_3,v_4$ are $R$-relaxed.
So by Lemma \ref{3_vertex_relaxed}, $v_5$ is a $(\geq 4)$-vertex.

By Claim 1, for every $i \in \{2,4,5\}$, either $f_i$ is not a 3-face, or $f_i$ is a 3-face but not receiving charge from $f$ due to (R2), (R3), or (R4).
Recall that $v_5$ is a $(\geq 4)$-vertex.
So $f$ only sends charge to $f_1$ and $f_3$.
Hence both $f_1$ and $f_3$ receive charge from $f$ due to (R2), for otherwise $\ch'(f) > (5-4) - 2 \cdot \frac{1}{2}=0$.

We shall show that Statement 3 holds for $i=1$.
Clearly, Statements 3(a), 3(b), 3(c) and 3(f) hold.
Since $f_1$ and $f_3$ receive charge from $f$ due to (R2), and $v_3$ and $v_2$ are $R$-relaxed, we know that Statements 3(d) and 3(e) hold by (R2). 
$\Box$

\medskip

\noindent{\bf Claim 3:} No 3-face receives charge from $f$ due to (R2), (R3), or (R4).

\noindent{\bf Proof of Claim 3:}
Suppose to the contrary that some 3-face, say $f_4$, receives charge from $f$ due to (R2), (R3), or (R4).
By Claim 2, $f_4$ is the only 3-face incident with $f$ receiving charge from $f$ due to (R2), (R3), or (R4).

We first suppose that $f_4$ receives charge from $f$ due to (R2).
Then $v_1$ or $v_3$ is non-$R$-relaxed.
By symmetry, we may assume that $v_1$ is non-$R$-relaxed.
So $v_2$ and $v_3$ are the only possible 3-vertices incident with $f$ by Lemma \ref{3_not_3}.
By Lemma \ref{3_vertex_relaxed}, at most one of $v_2$ and $v_3$ can be a 3-vertex.
So $\ch'(f) \geq (5-4) - \frac{1}{2} - \frac{1}{2} =0$, and the equality holds when $v_2$ or $v_3$ is a 3-vertex.
Since $f_4$ is a 3-face, $f$ is not a counterexample, a contradiction.

So $f_4$ receives charge from $f$ due to (R3) or (R4).
Hence $f$ sends at most $\frac{1}{3}$ to $f_4$.
So $f$ is incident with at least two 3-vertices, for otherwise $\ch'(f) \geq (5-4) - \frac{1}{3} - \frac{1}{2}>0$.
By Lemma \ref{3_vertex_relaxed} and \ref{3_cycle_relaxed}, $v_1$ and $v_3$ are 3-vertices, and $v_2$ is non-$R$-relaxed.
In particular, $\deg(v_2) \geq 4$, and if $\deg(v_2) \geq 5$, then $\deg(v_2) \geq 6$.

Now we suppose $\deg(v_2) \geq 5$.
So $\deg(v_2) \geq 6$.
Since $v_1$ and $v_3$ are 3-vertices, $f_1$ and $f_2$ are not 3-faces by Lemma \ref{3_not_3}.
So $v_2$ sends $\frac{1}{3}$ to $f$ by (R8).
Hence $\ch'(f) \geq (5-4) - \frac{1}{3} - 2 \cdot \frac{1}{2} + \frac{1}{3} \geq 0$.
Since $v_1$ and $v_3$ are two 3-vertices incident with $f$, $f$ is not a counterexample, a contradiction.

So $v_2$ is a 4-vertex.

Since $v_1$ and $v_3$ are 3-vertices and $v_4$ and $v_5$ are $R$-relaxed, Lemma \ref{3_vertex_relaxed} implies that for every $u \in (N(v_1) \cup N(v_3))-\{v_1,v_3,v_4,v_5\}$, $u$ is non-$R$-relaxed, so all edges incident with $u$ are not in $R$; moreover, since $v_2$ is a 4-vertex, Lemma \ref{relaxed_433} implies that for every neighbor $u$ of $v_2$ not incident with $f$, all edges incident with $u$ are not in $R$.
That is, for every $u \in (N(v_1) \cup N(v_2) \cup N(v_3))-\{v_1,v_3,v_4,v_5\}$, all edges incident with $u$ are not in $R$.

Hence, if $f_1$ is a 4-face, then the cycle bounding $f_1$ has $R$-length four, a contradiction; if $f_1$ is a 5-face, then the cycle bounding $f_1$ has $R$-length 5.
Since $v_1$ is a 3-vertex, $f_1$ is a $(\geq 5)$-face by Lemma \ref{3_not_3}.
Similarly, $f_2$ is a $(\geq 5)$-face, and if $f_2$ is a 5-face, then the cycle bounding $f_2$ has $R$-length 5.

Let $v_2'$ be the neighbor of $v_2$ incident with $f_1$ but not incident with $f$, and let $v_2''$ be the neighbor of $v_2$ incident with $f_2$ but not incident with $f$ such that $N_G(v_2)=\{v_1,v_3,v_2',v_2''\}$.

Recall that for every $u \in (N(v_1) \cup N(v_2) \cup N(v_3))-\{v_1,v_3,v_4,v_5\}$, all edges incident with $u$ are not in $R$.
So if $f_1$ is a $5$-face, then since there do not exist two cycles of $R$-length 5 sharing an edge, the face $f'$ sharing $v_2v_2'$ with $f_1$ is a $(\geq 6)$-face and sends $\frac{1}{4}$ unit of charge to $f$ via $v_2v_2'$ due to (R6).
If $f_1$ is a $(\geq 6)$-face, then $f_1$ sends $\frac{1}{4}$ unit of charge to $f$ via $v_1v_2$ due to (R5).
So either $f$ receives $\frac{1}{4}$ unit of charge from $f_1$ via $v_1v_2$ or from $f'$ via $v_2v_2'$.
Similarly, $f$ receives $\frac{1}{4}$ unit of charge from $f_2$ via $v_2v_3$ or from $f'$ via $v_2v_2''$.
Recall that $f$ sends at most $\frac{1}{3}$ unit of charge to $f_4$.
So $\ch'(f) \geq (5-4) - \frac{1}{3} - 2 \cdot \frac{1}{2} + 2 \cdot \frac{1}{4} > 0$, a contradiction.
$\Box$

\medskip

By Claim 3, $f$ only sends charge due to (R1).
If $f$ is incident with at most two 3-vertices, then $\ch'(f) \geq (5-4) - 2 \cdot \frac{1}{2}=0$, and the equality holds only when $f$ is incident with two 3-vertices, contradicting that $f$ is a counterexample.

So $f$ is incident with at least three 3-vertices.

Since no 3-vertex can be adjacent to two 3-vertices by Lemma \ref{3_vertex_relaxed}, we may without loss of generality assume that $v_1$, $v_3$ and $v_4$ are 3-vertices.
Moreover, by Lemma \ref{3_vertex_relaxed}, $v_3v_4$ is the only possible edge of $C$ belonging to $R$.
Since the $R$-length of $C$ is not 6, $v_3v_4$ is not in $R$.
Hence all edges of $C$ are not in $R$, and the $R$-length of $C$ equals 5.

Since $v_3$ and $v_4$ are 3-vertices, Lemma \ref{3_vertex_relaxed} implies that $f_2$ and $f_4$ are $(\geq 4)$-faces, and if $f_2$ or $f_4$ is a 4-face, then the cycle bounding this face must contain at most one edge in $R$ and hence has $R$-length 5, leading to a contradiction since the $R$-length of $C$ equals 5.
So $f_2$ and $f_4$ are $(\geq 5)$-faces.
If $f_2$ is a 5-face, then since $v_1$, $v_3$ and $v_4$ are 3-vertices, Lemma \ref{3_not_3} implies that $v_2v_3$ is the unique edge shared by $f_2$ and $f$, so Lemma \ref{5_faces_share_3_neighbor} (taking $F_1=f_2$ and $F_2=f$) implies that $\deg(v_1) \geq 4$ or $\deg(v_4) \geq 4$, a contradiction.
So $f_2$ is a $(\geq 6)$-face.
Similarly $f_4$ is a $(\geq 6)$-face.
So if $\deg(v_2)=4$, then $f_2$ sends $\frac{1}{4}$ unit of charge to $f$ via $v_2v_3$ due to (R5).
If $\deg(v_2) \geq 5$, then since $v_3$ and $v_4$ are 3-vertices, we know $\deg(v_3) \geq 6$ (by Lemma \ref{3_vertex_relaxed}) and $v_2$ sends $\frac{1}{3}$ unit of charge to $f$ due to (R8).
That is, $f$ receive at least $\frac{1}{4}$ unit of charge from $f_2$ via $v_2v_3$ or from $v_2$.
Similarly, $f$ receive at least $\frac{1}{4}$ unit of charge from $f_4$ via $v_4v_5$ or from $v_5$.
Hence $\ch'(f) \geq (5-4) - 3 \cdot \frac{1}{2} + 2 \cdot \frac{1}{4} \geq 0$.
Since $f$ is incident with at least two 3-vertices, $f$ is not a counterexample, a contradiction.
This proves the lemma.
\end{pf}

\bigskip

Let $F(G)$ be the set of faces of $G$.

\begin{lemma} \label{charge_0}
$\sum_{x \in V(G) \cup F(G)} \ch'(x) \geq 0$, and the equality holds only when $\ch'(x)=0$ for every $x \in V(G) \cup F(G)$.
\end{lemma}

\begin{pf}
This lemma immediately follows from Lemmas \ref{discharging_vertex}, \ref{discharging_3_face}, \ref{discharging_4_face}, \ref{discharging_6_face} and \ref{discharging_5_face}.
\end{pf}

\begin{lemma} \label{tight_no_6_face}
If $\sum_{x \in V(G) \cup F(G)} \ch'(x) = 0$, then every vertex has degree 3, 4, 5 or 6, and every face is a $3$-face, $4$-face or $5$-face.
\end{lemma}

\begin{pf}
By Lemma \ref{charge_0}, $\ch'(x)=0$ for every $x \in V(G) \cup F(G)$. 
By Lemma \ref{discharging_vertex}, every vertex has degree 3, 4, 5, or 6.
By Lemma \ref{discharging_6_face}, every face is a 3-face, 4-face, 5-face or 6-face.

\medskip

\noindent{\bf Claim 1:} Every 6-face is bounded by a cycle such that the degree of the vertices are 3, 3, 4, 3, 3, 4, listed in a cyclic order of the cycle bounding this face.

\noindent{\bf Proof of Claim 1:}
By Lemma \ref{discharging_vertex}, no vertex incident with a 6-face has degree at least six.
By Lemma \ref{discharging_6_face}, every 6-face is bounded by a cycle such that the degree of the vertices are 3, 3, 4, 3, 3, 4, listed in a cyclic order of the cycle bounding this face.
$\Box$

\medskip

To prove this lemma, it suffices to show that there is no 6-face.

Suppose to the contrary that there exists a 6-face $f$.
By Claim 1, $f$ is bounded by a cycle $v_1v_2v_3v_4v_5v_6v_1$, and the degree of the vertices are 3, 3, 4, 3, 3, 4, listed in a cyclic order of this cycle.

By symmetry, we may assume $\deg(v_1)=\deg(v_2)=3$ and $\deg(v_3)=4$.
Let $v_1'$ and $v_2'$ be the neighbors of $v_1$ and $v_2$ not incident with $f$, respectively.
Let $f_1$ be the face adjacent to $f$ sharing $v_1v_2$ with $f$.
Let $f_2$ be the face adjacent to $f$ sharing $v_2v_3$ with $f$.

We first suppose that $f_1$ is not a 6-face.
Since $v_1$ is a 3-vertex, $f_1$ is a 4-face or a 5-face by Lemma \ref{3_not_3}.
By Lemma \ref{3_vertex_relaxed}, $v_1'$ and $v_2'$ are non-$R$-relaxed vertices, so $v_1v_2$ is the only possible edge incident with $f_1$ belonging to $R$.
Hence the $R$-length of $f_1$ equals 5.
Since $v_2$ is a 3-vertex, $f_2$ is a $(\geq 4)$-face.
If $f_2$ is a 6-face, then since $\deg(v_2)=3$ and $\deg(v_3)=4$, Claim 1 implies that $v_2'$ is a 3-vertex, contradicting that $v_2'$ is non-$R$-relaxed.
So $f_2$ is not a 6-face.
Hence $f_2$ is a 4-face or a 5-face.
By Lemma \ref{relaxed_433}, the neighbor of $v_3$ incident with $f_2$ is non-$R$-relaxed.
Since $v_2'$ is also incident with $f_2$ and is non-$R$-relaxed, the $R$-length of $f_2$ equals 5.
Hence $f_1$ and $f_2$ are two adjacent faces of $R$-length 5.
Since there is no cycle of $R$-length 4 or 6, $f_1$ and $f_2$ share exactly one edge, a contradiction.

So $f_1$ is a 6-face.
By Claim 1, $v_1'$ and $v_2'$ are 4-vertices.
By Lemmas \ref{3_vertex_relaxed}, $v_2'$ and $v_3$ are non-$R$-relaxed.
So $v_2v_3 \not \in R$ and $v_2v_2' \not \in R$.
By Lemma \ref{relaxed_433}, the neighbor $v_3'$ of $v_3$ incident with $f_2$ is non-$R$-relaxed.
So $f_2$ is not a 3-face (by Lemma \ref{3_not_3}) or a 4-face (for otherwise the cycle bounding $f_2$ has $R$-length 4).
Since $v_2'$ and $v_3$ are 4-vertices, $f_2$ is not a 6-face by Claim 1.
So $f_2$ is a 5-face.

By Lemma \ref{relaxed_433} and Claim 1, the neighbor of $v_2'$ incident with $f_2$ other than $v_2$ is non-$R$-relaxed.
Recall that $v_3'$ is non-$R$-relaxed and $v_2v_3 \not \in R$ and $v_2v_2' \not \in R$.
So $v_2$ is the unique non-$R$-relaxed vertex incident with the 5-face $f_2$, and the cycle bounding $f_2$ has no edge in $R$. 
The former implies that $f_2$ is incident with at most one 3-vertex and no 3-face is adjacent to $f_2$ by Lemma \ref{3_cycle_relaxed}.
This contradicts Lemma \ref{discharging_5_face}.
\end{pf}

\begin{lemma} \label{tight_no_3_v}
If $\sum_{x \in V(G) \cup F(G)} \ch'(x) = 0$, then there is no 3-vertex.
\end{lemma}

\begin{pf}
Suppose to the contrary that there exists a 3-vertex $v$.
By Lemma \ref{discharging_vertex}, $v$ is incident with exactly two $(\geq 5)$-faces.
By Lemma \ref{tight_no_6_face}, these two $(\geq 5)$-faces are 5-faces.
By Lemma \ref{no_3_v_554_f}, $v$ is not incident with a 4-face.
So $v$ is incident with a 3-face, contradicting Lemma \ref{3_not_3}.
\end{pf}

\begin{lemma} \label{tight_no_5_face}
If $\sum_{x \in V(G) \cup F(G)} \ch'(x) = 0$, then there is no 5-face.
\end{lemma}

\begin{pf}
Suppose to the contrary that there exists a 5-face $f$.
By Lemma \ref{charge_0}, $\ch'(f)=0$.
By Lemma \ref{discharging_5_face}, $f$ is bounded by a cycle $C$.
Denote $C$ by $v_1v_2...v_5v_1$.
For every $i \in [5]$, let $f_i$ be the face sharing $v_iv_{i+1}$ with $f$, where $v_6=v_1$.
By Lemma \ref{tight_no_3_v}, every vertex in $C$ is a $(\geq 4)$-vertex.
So by Lemma \ref{discharging_5_face}, we may assume without loss of generality that 
    \begin{itemize}
        \item $f_1$ and $f_3$ are 3-faces and receive charge from $f$ due to (R2), (so $v_1,v_2,v_3,v_4$ are 4-vertices),
        \item the neighbor $v_2'$ of $v_2$ not incident with $f$ or $f_1$ is non-$R$-relaxed, and
        \item the neighbor $v_3'$ of $v_3$ not incident with $f$ or $f_3$ is non-$R$-relaxed,
        \item $v_5$ is non-$R$-relaxed.
    \end{itemize}
Since $v_2$ is a 4-vertex and $f_1$ is a 3-face, we know that $f_2$ is not a 3-face by Lemma \ref{relaxed_4_neighbors}.
By Lemma \ref{tight_no_6_face}, $f_2$ is a 4-face or a 5-face.

We first suppose that $f_2$ is a 5-face.
By Lemmas \ref{charge_0}, \ref{discharging_5_face} and \ref{tight_no_3_v}, we know that $f_2$ is incident with two 3-faces.
But it is impossible by Lemma \ref{3_not_3} since $v_2'$ and $v_3'$ are non-$R$-relaxed vertices incident with $f_2$, a contradiction.

So $f_2$ is not a 5-face.
By Lemma \ref{tight_no_6_face}, $f_2$ is a 3-face or a 4-face.
Since $v_2'$ is non-$R$-relaxed and incident with $f_2$, we know that $f_2$ is a 4-face.

Let $z_1$ be the vertex incident with $f_1$ other than $v_1$ and $v_2$.
Let $f_{12}$ be the face sharing $z_1v_2$ with $f_1$.

If $f_{12}$ is a 4-face, then since $f_1$ is a 3-face sharing an edge with $f_{12}$, Lemma \ref{3_4_face_relaxed} implies that all vertices incident with $f_{12}$ are $R$-relaxed; but $v_2'$ is non-$R$-relaxed and is incident with $f_{12}$, a contradiction.

So $f_{12}$ is not a 4-face.
By Lemmas \ref{no_3_cycle_share} and \ref{tight_no_6_face}, $f_{12}$ is a 5-face.
By Lemmas \ref{discharging_5_face} and \ref{tight_no_3_v}, $f_{12}$ is adjacent to two 3-faces with disjoint boundary receiving charge from $f_{12}$ due to (R2).
Since $v_2'$ is non-$R$-relaxed, $f_1$ is one of those two 3-faces.
Hence $z_1$ is a 4-vertex, and the neighbor $z^*$ of $z_1$ not incident with $f_{12} \cup f_1$ is non-$R$-relaxed.

Let $f'$ be the face sharing $z_1v_1$ with $f_1$.
Note that the three vertices $v_1,v_2,z_1$ incident with $f_1$ are 4-vertices.
So $f_1$ is adjacent to exactly one 4-face by Lemma \ref{discharging_3_face}.
Since $f$ and $f_{12}$ are 5-faces, $f'$ is a 4-face. 
Note that $z^*$ is incident with $f'$.
Since $f'$ is a 4-face and $f_1$ is a 3-face, $z^*$ is $R$-relaxed by Lemma \ref{3_4_face_relaxed}, a contradiction.
This proves the lemma.
\end{pf}

\begin{lemma} \label{tight}
If $\sum_{x \in V(G) \cup F(G)} \ch'(x) = 0$, then every face is a 4-face, and every vertex is a 4-vertex.
\end{lemma}

\begin{pf}
Let $f$ be a face.
By Lemmas \ref{tight_no_6_face} and \ref{tight_no_5_face}, $f$ is a 3-face or a 4-face.
If $f$ is a 3-face, then Lemmas \ref{no_3_cycle_share} and \ref{discharging_3_face} imply that $f$ is adjacent to a $(\geq 5)$-face, contradicting Lemmas \ref{tight_no_6_face} and \ref{tight_no_5_face}.
So $f$ is a 4-face.

This shows that every face is a 4-face.
Then every vertex is a 4-vertex by Lemma \ref{discharging_vertex}.
\end{pf}

\begin{lemma} \label{positive_charge}
$\sum_{x \in V(G) \cup F(G)} \ch'(x) > 0$.
\end{lemma}

\begin{pf}
Suppose to the contrary that $\sum_{x \in V(G) \cup F(G)} \ch'(x) \leq 0$.
By Lemma \ref{charge_0}, \linebreak $\sum_{x \in V(G) \cup F(G)} \ch'(x) = 0$.
By Lemma \ref{tight}, every vertex is a 4-vertex, and every face is a 4-face.
But it contradicts Lemma \ref{no_4_v_4_f}.
\end{pf}

\bigskip

Lemma \ref{positive_charge} implies that $\sum_{x \in V(G) \cup F(G)}\ch(x) = \sum_{x \in V(G) \cup F(G)}\ch'(x)>0$.
However, since $G$ is 2-cell embedded in a surface of Euler genus at most 2, $\sum_{x \in V(G) \cup F(G)}\ch(x) = \sum_{v \in V(G)} (\deg(v)-4) + \sum_{f \in F(G)}(\leng(f)-4) = -4(|V(G)|-|E(G)|+|F(G)|) \leq 0$ by Euler's formula, a contradiction.
This proves Theorem \ref{main_real}.

\bigskip

\noindent{\bf Acknowledgement:}
This work was conducted as part of the High School Research Program ``PReMa'' (Program for Research in Mathematics) at Texas A\&M University.
We thank MAA for supporting PReMa via the Dolciani Mathematics Enrichment Grant, and thank Kun Wang, the director of PReMa, for her organization.

\end{document}